\documentclass[a4paper,notitlepage,10pt]{article}%
\usepackage{amssymb}
\usepackage{graphicx}
\usepackage{amsmath}
\usepackage{amsthm}
\usepackage{amsfonts}
\usepackage{subfigure}%
\providecommand{\U}[1]{\protect\rule{.1in}{.1in}}

\theoremstyle{plain}
\newtheorem{theorem}{Theorem}[section]

\newtheorem{corollary}[theorem]{Corollary}

\newtheorem{lemma}[theorem]{Lemma}

\newtheorem{proposition}[theorem]{Proposition}

\theoremstyle{definition}

\newtheorem{remark}[theorem]{Remark}
\numberwithin{equation}{section}
\numberwithin{theorem}{section}

\let\footnote=\endnote

\ifx\pdfoutput\relax\let\pdfoutput=\undefined\fi
\newcount\msipdfoutput
\ifx\pdfoutput\undefined\else
\ifcase\pdfoutput\else
\ifx\paperwidth\undefined\else
\ifdim\paperheight=0pt\relax\else\pdfpageheight\paperheight\fi
\ifdim\paperwidth=0pt\relax\else\pdfpagewidth\paperwidth\fi
\fi\fi\fi
\begin{document}

\title{A curve of positive solutions for an indefinite sublinear Dirichlet problem}
\author{Uriel Kaufmann\thanks{FaMAF, Universidad Nacional de C\'{o}rdoba, (5000)
C\'{o}rdoba, Argentina. \textit{E-mail address: }kaufmann@mate.uncor.edu} ,
Humberto Ramos Quoirin\thanks{Universidad de Santiago de Chile, Casilla 307,
Correo 2, Santiago, Chile. \textit{E-mail address: }humberto.ramos@usach.cl} ,
Kenichiro Umezu\thanks{Department of Mathematics, Faculty of Education,
Ibaraki University, Mito 310-8512, Japan. \textit{E-mail address:
}kenichiro.umezu.math@vc.ibaraki.ac.jp}
\and \noindent}
\maketitle

\begin{abstract}
We investigate the existence of a curve $q\mapsto u_{q}$, with $q\in(0,1)$, of
positive solutions for the problem
\[%
\begin{cases}
-\Delta u=a(x)u^{q} & \mbox{ in }\Omega,\\
u=0 & \mbox{ on }\partial\Omega,
\end{cases}
\leqno{(P_{a,q})}
\]
where $\Omega$ is a bounded and smooth domain of $\mathbb{R}^{N}$ and
$a:\Omega\rightarrow\mathbb{R}$ is a sign-changing function (in which case the
strong maximum principle does not hold). In addition, we analyze the
asymptotic behavior of $u_{q}$ as $q\rightarrow0^{+}$ and $q\rightarrow1^{-}$.
We also show that in some cases $u_{q}$ is the ground state solution of
$(P_{a,q})$.
As a byproduct, we obtain existence results for a singular and indefinite
Dirichlet problem. Our results are mainly based on bifurcation and
sub-supersolutions methods.

\end{abstract}


\section{Introduction}

Let $\Omega$ be a bounded and smooth domain of $\mathbb{R}^{N}$ with $N\geq1$,
and $0<q<1$. In this article we proceed with the investigation of the problem
\[%
\begin{cases}
-\Delta u=a(x)u^{q} & \mbox{ in }\Omega,\\
u=0 & \mbox{ on }\partial\Omega,
\end{cases}
\leqno{(P_{a,q})}
\]
where $\Delta$ is the usual Laplacian in $\mathbb{R}^{N}$. Throughout this
paper, unless otherwise stated, we assume that $r>N$ and $a\in L^{r}\left(
\Omega\right)  $ is such that $\Omega_{+}(a)\neq\emptyset$ and
$|(\mathrm{supp}\,a^{+})\setminus\Omega_{+}(a)|=0$. Here $\Omega_{+}(a)$ is
the largest open subset of $\Omega$ where $a>0$ a.e., $a^{\pm}:=\max(\pm
a,0)$, and $|A|$ stands for the Lebesgue measure of $A\subset\mathbb{R}^{N}$.
If there is no ambiguity we simply write $\Omega_{+}$ instead of $\Omega
_{+}(a)$.

Let us set
\[
W_{D}^{2,r}(\Omega):=\{u\in W^{2,r}(\Omega):u=0\text{ on }\partial\Omega\}.
\]
By a \textit{nonnegative} \textit{solution} of $(P_{a,q})$ we mean a function
$u\in W_{D}^{2,r}\left(  \Omega\right)  $ (and thus $u\in C^{1}(\overline
{\Omega})$) that satisfies the equation for the weak derivatives and $u\geq0$
in $\Omega$. If, in addition, $u>0$ in $\Omega$, then we call it a
\textit{positive solution} of $(P_{a,q})$.

We are mostly interested in the case where $a$ changes sign, i.e. $|\text{supp
}a^{\pm}|>0$. This case turns out to be the most interesting and challenging
one, since neither the strong maximum principle nor Hopf's Lemma apply and
consequently nonnegative solutions of $(P_{a,q})$ do not necessarily belong
to
\[
\mathcal{P}^{\circ}:=\left\{  u\in C_{0}^{1}(\overline{\Omega}):u>0\text{ in
}\Omega\text{ \ and }\frac{\partial u}{\partial\nu}<0\text{ on }\partial
\Omega\right\}  ,
\]
the interior of the positive cone of
\[
C_{0}^{1}(\overline{\Omega}):=\left\{  u\in C^{1} (\overline{\Omega
}):u=0\text{ on }\partial\Omega\right\}  .
\]
In fact, a nonnegative solution of $(P_{a,q})$ may vanish in parts of $\Omega$
and its normal derivative may vanish in parts of $\partial\Omega$. This
phenomenon provides a rich structure to the nonnegative solutions set of
$(P_{a,q})$. It is our purpose in this article to better understand this structure.

Let $\phi\in W_{D}^{2,r}\left(  \Omega\right)  $ be the unique solution of%
\[%
\begin{cases}
-\Delta\phi=a(x) & \mbox{ in }\Omega,\\
\phi=0 & \mbox{ on }\partial\Omega,
\end{cases}
\leqno{(P_a)}
\]
and $\mathcal{S}:L^{r}\left(  \Omega\right)  \rightarrow W_{D}^{2,r}\left(
\Omega\right)  $ be the corresponding solution operator, i.e. $\mathcal{S}%
(a)=\phi$. If $\phi>0$ in $\Omega$ then $\phi^{0}=1$ is well defined so that,
in this case, we set $(P_{a,0}):=(P_{a})$.

Let us mention that the existence of nonnegative solutions of $(P_{a,q})$ was
considered in detail in \cite{BPT} (see also \cite{PT}), assuming that $a$ is
H\"{o}lder continuous. In particular, it was shown in \cite{BPT} that
$(P_{a,q})$ admits a nontrivial nonnegative solution for all $q\in\left(
0,1\right)  $, and that uniqueness of nontrivial nonnegative solutions does
\textit{not} hold in general for $(P_{a,q})$.

However, the existence of positive solutions for $(P_{a,q})$ is a more
delicate issue which has been addressed by very few papers. It was first
proved in \cite[Theorem 4.4]{jesusultimo} (see also \cite[Theorem 3.7]{royal})
that if $\mathcal{S}(a)\in\mathcal{P}^{\circ}$ then $(P_{a,q})$ has a positive
solution (which may \textit{not} belong to $\mathcal{P}^{\circ}$). We note
that this condition is not sharp, since there exists $a$ such that
$\mathcal{S}(a)<0$ in $\Omega$ but $(P_{a,q})$ has a positive solution for
some $q\in(0,1)$ (see \cite[Section 1]{nodea}). Later on, in \cite{nodea}, the
authors studied $(P_{a,q})$ in the one-dimensional and radial cases,
establishing several sufficient conditions on $a$ (as well as some necessary
ones) for the existence of a positive solution of $(P_{a,q})$. Some of these
results were then extended to the case of a smooth bounded domain in
\cite{ans}.

More recently, we have proved in \cite{krqu} that if $\Omega_{+}$ has finitely
many connected components and $q$ is close enough to $1$, then any nontrivial
nonnegative solution of $(P_{a,q})$ belongs to $\mathcal{P}^{\circ}$, so that
in this situation
$(P_{a,q})$ has \textit{exactly} one nontrivial nonnegative solution, which
\textit{in addition} belongs to $\mathcal{P}^{\circ}$ (see \cite[Corollary
1.5]{krqu}). This uniqueness and positivity result was proved via a continuity
argument inspired by \cite[Theorem 4.1]{J} (see also \cite{K}), which is based
on the fact that the strong maximum principle applies to $(P_{a,q})$ if $q=1$.
Furthermore, it also applies to the Neumann counterpart of $(P_{a,q})$, cf.
\cite[Theorem 1.7 and Corollary 1.8]{krqu}. We refer to \cite{BPT2} and our
recent work \cite{neumann} for a detailed study of $(P_{a,q})$ under Neumann
boundary conditions.
Let us mention that similar problems in the case $\Omega=\mathbb{R}^{N}$ have
been considered in \cite{AL09, BK92,Sp83}.

Regarding other uniqueness results for $(P_{a,q})$, let us recall the following:

\setcounter{theorem}{-1}

\begin{theorem}
\label{t0} \strut

\begin{enumerate}
\item \textrm{\cite[Theorem 2.1]{DS} $(P_{a,q})$ has at most one positive
solution. }

\item \textrm{\cite[Theorem 2.1]{BPT} Assume $a\in C^{\alpha}(\overline
{\Omega})$, $\alpha\in\left(  0,1\right)  $. If $\Omega_{+}$ is smooth and has
finitely many connected components, then $(P_{a,q})$ has at most one solution
positive in $\Omega_{+}$. }
\end{enumerate}
\end{theorem}

\begin{remark}
\label{r0}From the proof of \cite[Theorem 2.1]{BPT} one may easily check that
this result remains valid for $a\in L^{r}\left(  \Omega\right)  $, $r>N$,
under the additional condition $(H_{+})$ below.
\end{remark}

In view of the above theorem, whenever $(P_{a,q})$ has a positive solution, we
denote it by $u(q)$.
Let
\[
\mathcal{I}_{a}:=\{q\in(0,1):(P_{a,q})\text{ has a solution }u\in
\mathcal{P}^{\circ}\},
\]
i.e.
\[
\mathcal{I}_{a}=\{q\in(0,1):u(q)\in\mathcal{P}^{\circ}\}.
\]
Our results shall be established under different conditions (most of them
technical), which are listed below:
\[
\left\{
\begin{array}
[c]{l}%
\Omega_{+}\text{ consists of \textit{finitely many} connected components and
}\\
\partial\Omega_{+}\text{ satisfies an inner sphere condition with respect to
}\Omega_{+},
\end{array}
\right.  \leqno{(H_{+})}
\]%
\[
\left\{
\begin{array}
[c]{l}%
\Omega_{+}\text{ is connected and }\\
\partial\Omega_{+}\text{ satisfies an inner sphere condition with respect to
}\Omega_{+},
\end{array}
\right.  \leqno{(H_{+}')}
\]

\[
\left\vert a(x)\right\vert \leq Cd(x,\partial\Omega)^{\alpha}\text{
\ \negthinspace a.e.\ in }\Omega_{\rho_{0}},\text{ for some }\rho_{0}>0\text{
and }\alpha>1-\frac{1}{N},\leqno{(H_1)}
\]
where
\[
\Omega_{\rho}=:\{x\in\Omega:d(x,\partial\Omega)<\rho\}.
\]

Given a subdomain $\Omega^{\prime}\subseteq\Omega$ such that $a^{+}%
\not \equiv 0$ in $\Omega^{\prime}$, we denote by $\lambda_{1}(a,\Omega
^{\prime})$ (or simply $\lambda_{1}\left(  a\right)  $ when $\Omega^{\prime
}=\Omega$) the first positive eigenvalue of the problem
\[%
\begin{cases}
-\Delta\phi=\lambda a(x)\phi & \mbox{ in }\Omega^{\prime},\\
\phi=0 & \mbox{ on }\partial\Omega^{\prime},
\end{cases}
\]
and by $\phi_{1}=\phi_{1}(a,\Omega^{\prime})$ the positive eigenfunction
associated to $\lambda_{1}(a,\Omega^{\prime})$ with $\int_{\Omega^{\prime}%
}\phi_{1}^{2}=1$.

We are now in position to state our main results. First we study the
properties of nonnegative ground state solutions of $(P_{a,q})$, i.e.
nonnegative global minimizers
of $I_{q}:H_{0}^{1}(\Omega)\rightarrow\mathbb{R}$, where
\[
I_{q}(u):=\frac{1}{2}\int_{\Omega}|\nabla u|^{2}-\frac{1}{q+1}\int_{\Omega
}a|u|^{q+1}.
\]
Some properties of such solutions have been already proved in \cite[Theorem
2.3]{BPT}. We complement it in the following result:

\begin{theorem}
\label{mt2a} $(P_{a,q})$ has a unique nonnegative ground state solution
$U_{q}$ for every $q\in(0,1)$, which satisfies $U_{q}>0$ in $\Omega_{+}$ and
$q\mapsto U_{q}$ is continuous from $(0,1)$ to $W_{D}^{2,r}(\Omega)$. Moreover
$U_{q}=u(q)\in\mathcal{P}^{\circ}$ for $q<1$ sufficiently close to $1$, and
$U_{q}$ has the following asymptotic behavior as $q\rightarrow1^{-}$:

\begin{itemize}
\item $U_{q}\rightarrow0$ in$\ W_{D}^{2,r}(\Omega)$ if $\lambda_{1}(a)>1$.

\item $\Vert U_{q}\Vert_{C(\overline{\Omega})}\rightarrow\infty$ if
$\lambda_{1}(a)<1$.
\end{itemize}

\noindent In addition, the following assertions hold:

\begin{enumerate}
\item If $q_{n} \to0^{+}$ then, up to a subsequence, $U_{q_{n}} \to U_{0}$ in
$C_{0}^{1}(\overline{\Omega})$, where $U_{0}$ is a nonnegative global
minimizer of $I_{0}$. In particular, if $0\not \equiv \mathcal{S}(a)\geq0$ in
$\Omega$, then $U_{q}\rightarrow\mathcal{S}(a)$ in $C_{0}^{1}(\overline
{\Omega})$ as $q\rightarrow0^{+}$. If $\mathcal{S}(a)>0$ in $\Omega$, then
this convergence holds in $W_{D}^{2,r}(\Omega)$.

\item If $(H_{1})$ holds and $\lambda_{1}(a)=1$, then $U_{q}\rightarrow
t_{\ast}\phi_{1}$ in $W_{D}^{2,r}(\Omega)$ as $q\rightarrow1^{-}$, where
$\phi_{1}=\phi_{1}(a,\Omega)$ and
\[
t_{\ast}:=\exp\left[  -\frac{\int_{\Omega}a\left(  x\right)  \phi_{1}^{2}%
\log\phi_{1}}{\int_{\Omega}a\left(  x\right)  \phi_{1}^{2}}\right]  .
\]
More precisely, $U_{q}$ bifurcates from $(1,t_{\ast}\phi_{1})$ to the region
$q<1$, and moreover, $q\mapsto U_{q}$ is continuous from $(q_{1},1]$ to
$W_{D}^{2,r}(\Omega)$, for some $q_{1}\in(0,1)$, where we set $U_{1}:=t_{\ast
}\phi_{1}$.


\item If $(H_{+})$ holds then $U_{q}$ is the maximal nonnegative solution of
$(P_{a,q})$, i.e. $U_{q}\geq u$ for any nonnegative solution $u$ of
$(P_{a,q})$. Furthermore,%
\begin{equation}
\mathcal{I}_{a}=\{q\in(0,1):U_{q}\in\mathcal{P}^{\circ}\}, \label{iia}%
\end{equation}
and $\mathcal{I}_{a}$ is open.
\end{enumerate}
\end{theorem}

\begin{remark}
\label{remgs} \strut

\begin{enumerate}
\item Whenever
$U_{q}>0$ in $\Omega$, we have, by Theorem \ref{t0} (i), $U_{q}=u(q)$. This
equality also holds whenever $u\left(  q\right)  $ exists and $(H_{+})$ is
satisfied, as a consequence of the first assertion in Theorem \ref{mt2a} (iii).

\item We can give a better asymptotic estimate for
$U_{q}$ as $q\rightarrow1^{-}$ when $(H_{1})$ holds and $\lambda_{1}(a)\neq1$.
Indeed, in this case Theorem \ref{mt2a} (ii) and a rescaling argument yield
that
\[
U_{q}\sim\lambda_{1}(a)^{-\frac{1}{1-q}}\,t_{\ast}\phi_{1}\quad\mbox{as}\quad
q\rightarrow1^{-},
\]
i.e. $\lambda_{1}(a)^{\frac{1}{1-q}}U_{q}\rightarrow t_{\ast}\phi_{1}$ in
$W_{D}^{2,r}(\Omega)$ as $q\rightarrow1^{-}$, see Corollary \ref{cor:0to1}.
\end{enumerate}
\end{remark}

Our second result slightly improves \cite[Theorem 4.4]{jesusultimo} showing
that $\mathcal{S}(a)>0$ in $\Omega$ is enough to get the existence of $u(q)$
for all $q\in(0,1)$. Moreover, it provides the asymptotic behavior of $u(q)$
as $q\rightarrow0^{+}$, as well as sufficient conditions to have
$u(q)\in\mathcal{P}^{\circ}$ for \textit{every} $q\in(0,1)$:

\begin{theorem}
\label{mt2b} Assume $\mathcal{S}(a)>0$ in $\Omega$. Then, there exists a
(unique) positive solution $u(q)$ of $(P_{a,q})$ for every $q\in(0,1)$.
Moreover, $u(q)\rightarrow\mathcal{S}(a)$ in $W_{D}^{2,r}(\Omega)$ as
$q\rightarrow0^{+}$, and the following assertions hold:

\begin{enumerate}
\item If $(H_{+})$ holds then $U_{q}=u(q)$ for every $q\in(0,1)$, and
$q\mapsto u(q)$ is continuous from $[0,1)$ to $W_{D}^{2,r}(\Omega)$, where we
set $u(0):=\mathcal{S}(a)$.

\item If $\mathcal{S}(a)\in\mathcal{P}^{\circ}$, then $U_{q}=u(q)\in
\mathcal{P}^{\circ}$ for $q>0$ sufficiently small.

\item Assume one of the following conditions:

\begin{enumerate}
\item $a\geq0$ in $\Omega_{\rho_{0}}$ for some $\rho_{0}>0$,

\item $\Omega$ is a ball and $a$ is radial,

\item $\mathcal{S}\left(  a\right)  \in\mathcal{P}^{\circ}$ and $(H_{+}%
^{\prime})$ holds.
\end{enumerate}

Then $u(q)\in\mathcal{P}^{\circ}$ for all $q\in(0,1)$, i.e. $\mathcal{I}%
_{a}=\left(  0,1\right)  $.
\end{enumerate}
\end{theorem}

The implicit function theorem also provides us with the following result:

\begin{theorem}
\label{mt2c} If $(H_{1})$ holds then $\mathcal{I}_{a}$ is open, and $u(q)$ is
asymptotically stable for $q\in\mathcal{I}_{a}$.
\end{theorem}

\begin{figure}[tbh]
\label{f1}
\par
\begin{center}
\includegraphics[scale=0.25]{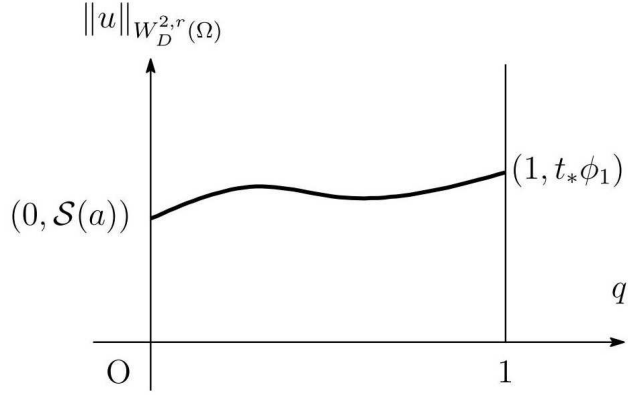}
\end{center}
\caption{The curve of positive solutions in the case $\lambda_{1}(a)=1$.}%
\label{fig17_0618b}%
\end{figure}

\begin{figure}[h]
\centerline{
\includegraphics[scale=0.2]{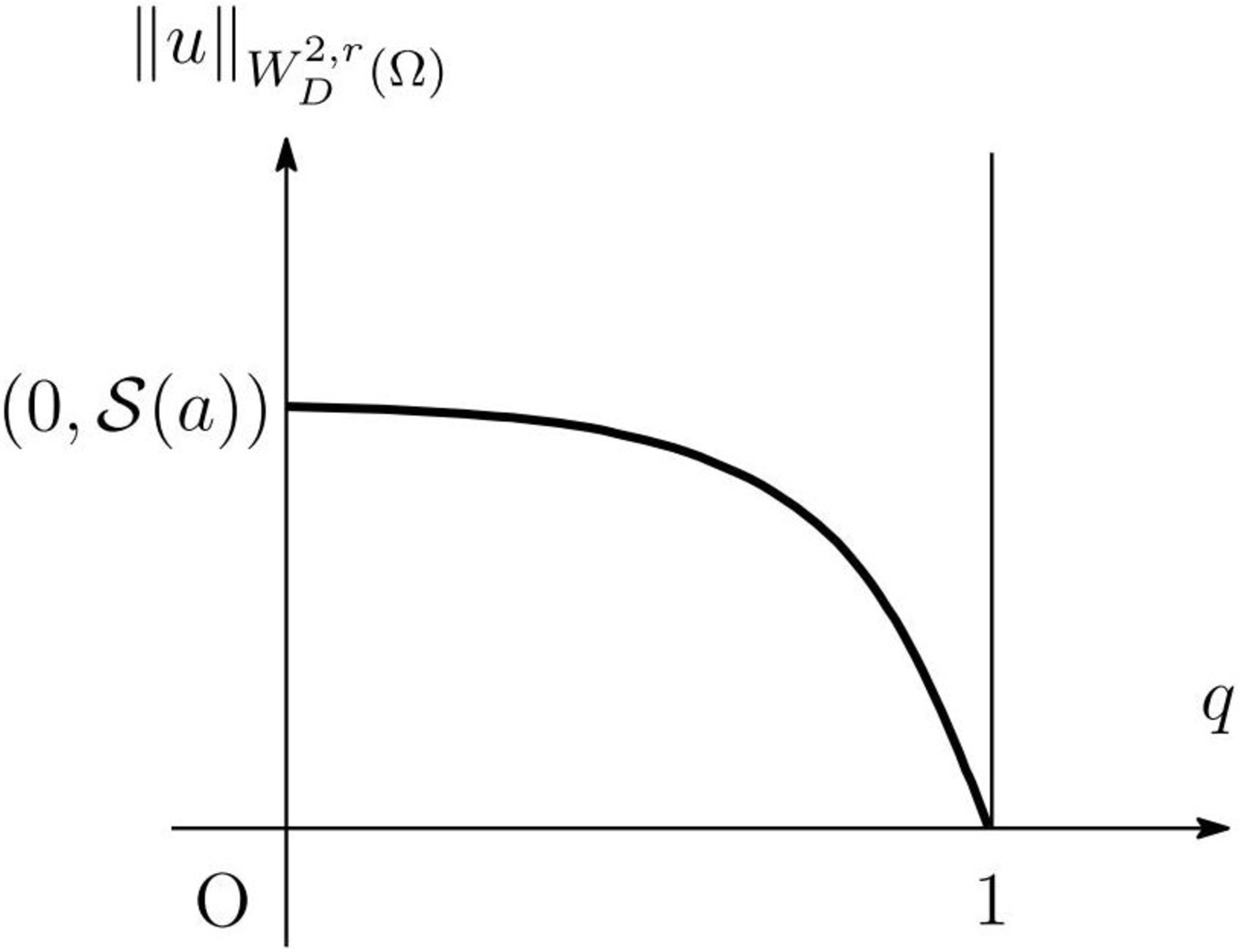} \hskip0.35cm
\includegraphics[scale=0.2] {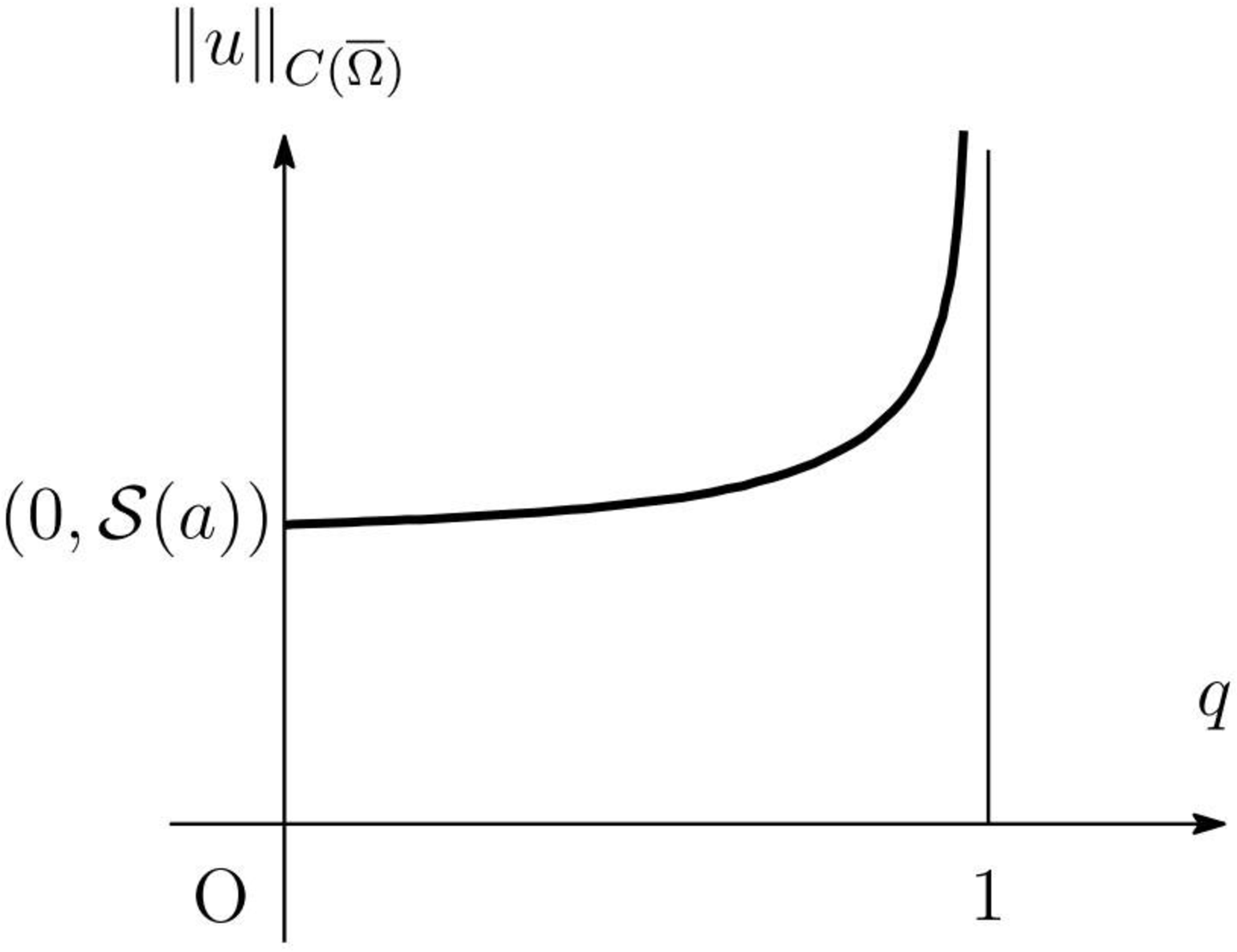}
} \centerline{(i) \hskip6.1cm (ii) }\caption{The curve of positive solutions
emanating from $(0,\mathcal{S}(a))$: (i) The case $\lambda_{1}(a)>1$. (ii) The
case $\lambda_{1}(a)<1$.}%
\label{fig:double}%
\end{figure}

As a direct consequence of the above theorems we obtain the next
result:

\begin{corollary}
\label{cormt2} Assume $(H_{+})$, $(H_{1})$ and $\mathcal{S}(a)>0$ in $\Omega$.

\begin{enumerate}
\item If $\lambda_{1}(a)=1$, then $U_{q}=u(q)>0$ in $\Omega$ for every
$q\in(0,1)$, and $q\mapsto U_{q}$ is continuous from $[0,1]$ to $W_{D}%
^{2,r}(\Omega)$, where we set $u(0):=\mathcal{S}\left(  a\right)  $ and
$u(1):=t_{\ast}\phi_{1}$ (see Figure \ref{fig17_0618b}). Moreover, if in
addition some of the conditions in Theorem \ref{mt2b} (iii) hold,
then $u\left(  q\right)  \in\mathcal{P}^{\circ}$ for all $q\in\left(
0,1\right)  $.

\item If $\lambda_{1}(a)\not =1$, then the curve of positive solutions in (i)
bifurcates at $q=1$ from zero or infinity, in accordance with the sign of
$\lambda_{1}(a)-1$ (see Figure \ref{fig:double}).
\end{enumerate}
\end{corollary}

\begin{remark}
\strut

\begin{enumerate}
\item From Proposition \ref{pgs2} below, we see that $(H_{1})$ can be removed
in Corollary \ref{cormt2} (ii).

\item We believe that the conclusions of Theorems \ref{mt2a} (ii) and
\ref{mt2c} remain valid even if $(H_{1})$ does not hold.

\item Let us point out that, to the best of our knowledge, the asymptotic
behavior of $u_{q}$ as $q\rightarrow0^{+}$ and $q\rightarrow1^{-}$, given by
Theorems \ref{mt2a} and \ref{mt2b}, are new even in the case $a\geq0$.

\item Theorem \ref{mt2b} (i) can be complemented as follows: Let $q_{0}%
\in(0,1)$. Then the set of nontrivial nonnegative solutions of $(P_{a,q})$ in
a $\mathbb{R}\times W_{D}^{2,r}(\Omega)$-neighborhood of $(q_{0},u(q_{0}))$ is
precisely given by $\{(q,u(q)):q\in(q_{0}-\delta_{0},q_{0}+\delta_{0})\}$, for
some $\delta_{0}>0$. Indeed, we first show that if $(q,u)$ is sufficiently
close to $(q_{0},u(q_{0}))$ and $u$ is a nontrivial nonnegative solution of
$(P_{a,q})$, then $u>0$ in $\Omega_{+}$. To this end, assume by contradiction
that $(q_{n},u_{n})\rightarrow(q_{0},u(q_{0}))$ in $\mathbb{R}\times
W_{D}^{2,r}(\Omega)$ and $u_{n}$ is a nontrivial nonnegative solution of
$(P_{a,q_{n}})$, but $u_{n}\left(  x_{n}\right)  =0$ for some $x_{n}\in
\Omega_{+}$. Since $(H_{+})$ holds, we may assume that $x_{n}\in\Omega
^{\prime}$, where $\Omega^{\prime}$ is a fixed connected component of
$\Omega_{+}$. From the strong maximum principle, $u_{n}\equiv0$ in
$\Omega^{\prime}$ for all $n$, and thus $u\left(  q_{0}\right)  \equiv0$ in
$\Omega^{\prime}$ which is not possible. Therefore, since $u\left(  q\right)
$ exists, Remark \ref{r0} ensures that $u=u(q)$. In relation with this result,
we note that if
in addition $(H_{+}^{\prime})$ holds, then $u(q)$ is the unique nontrivial
nonnegative solution of $(P_{a,q})$.
\end{enumerate}
\end{remark}

As a byproduct of the above theorems, we obtain an existence result for a
singular and indefinite Dirichlet problem. Although singular problems like
$\left(  P_{s}\right)  $ below have been extensively studied when $a\geq0$ in
$\Omega$ (see e.g. \cite{tartar,pino,royal,gomes,lair,lazer} and references
therein), as far as we know, the indefinite case has only been considered
recently in \cite{jmaa} (see also \cite{ana} for a similar problem with the
one-dimensional $p$-Laplacian). The following corollary complements (with a
different approach) some of the results presented there.

\begin{corollary}
\label{sing}Assume that $(H_{1})$ holds and $\mathcal{S}\left(  a\right)
\in\mathcal{P}^{\circ}$. Then there exists $\gamma_{0}>0$ such that the
singular problem%
\[%
\begin{cases}
-\Delta u=a(x)u^{-\gamma} & \mbox{ in }\Omega,\\
u=0 & \mbox{ on }\partial\Omega,
\end{cases}
\leqno{(P_{s})}
\]
has a solution $u\in\mathcal{P}^{\circ}$ for all $\gamma\in\left(
0,\gamma_{0}\right)  $. Moreover, there exists $s>N$ such that $\gamma\mapsto
u(\gamma)$ is continuous from $\left[  0,\gamma_{0}\right)  $ to $W_{D}%
^{2,s}(\Omega)$, where $u\left(  0\right)  :=\mathcal{S}\left(  a\right)  $.
\end{corollary}

\begin{remark}
It was proved in \cite[Corollary 4.6]{jmaa} that the condition $\mathcal{S}%
\left(  a\right)  >0$ in $\Omega$ is necessary for the existence of solutions
of $\left(  P_{s}\right)  $ lying in $W_{D}^{2,s}(\Omega)$ with $s>N$.
Therefore, we see that the sufficient condition imposed in Corollary
\ref{sing}, namely $\mathcal{S}\left(  a\right)  \in\mathcal{P}^{\circ}$, is
\textquotedblleft almost\textquotedblright\ sharp.\smallskip
\end{remark}

The outline of this article is the following: in Section \ref{sec:g} we obtain
some
properties of nonnegative ground state solutions of $(P_{a,q})$. Sections
\ref{sec:b1} and \ref{sec:b2} are devoted to the use of the implicit function
theorem and a bifurcation analysis of $(P_{a,q})$, where $q$ is regarded as a
bifurcation parameter. Finally, Section \ref{sec:p} provides some additional
results obtained by the sub-supersolutions method and the proofs of our main theorems.

\subsubsection*{Notation}

\begin{itemize}
\setlength{\itemsep}{0.0cm}

\item For any $f\in L^{1}(\Omega)$ the integral $\int_{\Omega}f$ is considered
with respect to the Lebesgue measure, whereas for any $g\in L^{1}%
(\partial\Omega)$ the integral $\int_{\partial\Omega}g$ is considered with
respect to the surface measure.

\item The weak convergence is denoted by $\rightharpoonup$.

\item The positive and negative parts of a function $u$ are defined by
$u^{\pm}:=\max\{\pm u,0\}$.

\item $|\cdot|$ stands for both the Lebesgue measure and the surface measure.

\item The characteristic function of a set $A\subset\mathbb{R}^{N}$ is denoted
by $\chi_{A}$.

\item $2^{*}=\frac{2N}{N-2}$ if $N \geq3$ and $2^{*}=\infty$ if $N=1,2$.

\end{itemize}

\section{The ground state solution}

\label{sec:g}

The following lemma will be frequently used in the sequel:

\begin{lemma}
\label{lem:c}\strut

\begin{enumerate}
\item Let $q_{n}\in(0,1)$ with $q_{n}\rightarrow q_{0}\in(0,1]$, and
$v(q_{n})$ be nontrivial nonnegative solutions of $(P_{a,q_{n}})$. If
$v(q_{n})\rightarrow v(q_{0})$ in $C(\overline{\Omega})$, then $v(q_{n}%
)\rightarrow v(q_{0})$ in $W_{D}^{2,r}(\Omega)$.

\item Let $q_{n}\rightarrow0^{+}$, and $v(q_{n})$ be nontrivial nonnegative
solutions of $(P_{a,q_{n}})$. If $v(q_{n})\rightarrow\mathcal{S}(a)$ in
$C(\overline{\Omega})$ and $\mathcal{S}(a)>0$ in $\Omega$, then $v(q_{n}%
)\rightarrow\mathcal{S}(a)$ in $W_{D}^{2,r}(\Omega)$.
\end{enumerate}
\end{lemma}

\quad

\textit{Proof.} We prove assertion (i). The proof of (ii) is carried out in a
similar way, so we omit it.

It is clear that
\[%
\begin{cases}
-\Delta(v(q_{n})-v(q_{0}))=a\left(  x\right)  (v(q_{n})^{q_{n}}-v(q_{0}%
)^{q_{0}}) & \mbox{ in }\Omega,\\
v(q_{n})-v(q_{0})=0 & \mbox{ on }\partial\Omega.
\end{cases}
\]
Using the Lebesgue dominated convergence theorem, we shall show that
\begin{equation}
\Vert a\left(  x\right)  (v(q_{n})^{q_{n}}-v(q_{0})^{q_{0}})\Vert
_{L^{r}(\Omega)}\rightarrow0. \label{W2rcoa+}%
\end{equation}
Once this is done, we have by elliptic regularity that $v(q_{n})\rightarrow
v(q_{0})$ in $W_{D}^{2,r}(\Omega)$ as $q_{n}\rightarrow q_{0}$, as desired.

Let us first deduce that
\begin{equation}
v(q_{n})^{q_{n}}\rightarrow v(q_{0})^{q_{0}}\quad\mbox{in}\quad\Omega.
\label{pwcq0}%
\end{equation}
Indeed, this is clear for $x\in\Omega$ such that $v(q_{0})(x)>0$. Now, if
$v(q_{0})(x)=0$, then we have two possibilities for $v(q_{n})$:
\[
v(q_{n})^{q_{n}}=\left\{
\begin{array}
[c]{ll}%
0 & \mbox{if}\quad v(q_{n})(x)=0,\\
e^{q_{n}\log v(q_{n})} & \mbox{if}\quad v(q_{n})(x)>0.
\end{array}
\right.
\]
Since $v(q_{n})(x)\rightarrow v(q_{0})(x)=0$ and $q_{0}>0$, it follows that
\[
v(q_{n})^{q_{n}}\rightarrow0=0^{q_{0}}=v(q_{0})^{q_{0}},
\]
and thus, \eqref{pwcq0} has been verified. This implies that
\begin{equation}
a\left(  x\right)  (v(q_{n})^{q_{n}}-v(q_{0})^{q_{0}})\rightarrow
0\quad\mbox{a.e. in}\quad\Omega. \label{W2rco+}%
\end{equation}
On the other side, since $v(q_{n})\rightarrow v(q_{0})$ in $C(\overline
{\Omega})$, we have $v(q_{n})\leq C_{1}$ on $\overline{\Omega}$ for some
$C_{1}>0$ independent of $n$. Hence, we infer that
\[
v(q_{n})^{q_{n}}=e^{q_{n}\log v(q_{n})}\leq e^{\log C_{1}}\quad\mbox{if}\quad
v(q_{n})(x)>0,
\]
and consequently, for some $C>0$
\begin{equation}
|a\left(  x\right)  (v(q_{n})^{q_{n}}-v(q_{0})^{q_{0}})|^{r}\leq
C|a(x)|^{r}\quad\mbox{a.e.\ in}\quad\Omega. \label{W2rco02+}%
\end{equation}
From \eqref{W2rco+} and \eqref{W2rco02+}, the Lebesgue dominated convergence
theorem yields \eqref{W2rcoa+}. $\blacksquare$ \newline

We recall that $I_{q}:H_{0}^{1}(\Omega)\rightarrow\mathbb{R}$ is given by
\[
I_{q}(u):=\frac{1}{2}\int_{\Omega}|\nabla u|^{2}-\frac{1}{q+1}\int_{\Omega
}a|u|^{q+1}%
\]
for $q\in\lbrack0,1)$.

\begin{proposition}
\label{pgs} $I_{q}$ has a unique nonnegative global minimizer $U_{q}$ for
every $q\in(0,1)$. In addition:

\begin{enumerate}
\item $U_{q}>0$ in $\Omega_{+}$ and $q\mapsto U_{q}$ is continuous from
$(0,1)$ to $W_{D}^{2,r}(\Omega)$.

\item There exists $q_{0}\in(0,1)$ such that $U_{q}\in\mathcal{P}^{\circ}$ for
$q\in(q_{0},1)$.

\item If $q_{n}\rightarrow0^{+}$ then, up to a subsequence, $U_{q_{n}%
}\rightarrow U_{0}$ in $C_{0}^{1}(\overline{\Omega})$, where $U_{0}$ is a
nonnegative global minimizer of $I_{0}$. In particular, if $0\not \equiv
\mathcal{S}(a)\geq0$ in $\Omega$, then $U_{q}\rightarrow\mathcal{S}(a)$ in
$C_{0}^{1}(\overline{\Omega})$ as $q\rightarrow0^{+}$. If $S(a)>0$ in $\Omega$
then this convergence holds in $W_{D}^{2,r}(\Omega)$.

\end{enumerate}
\end{proposition}

\textit{Proof}. By a standard minimization argument, one may easily prove the
existence of a global minimizer of $I_{q}$. Moreover, there is a $1$ to $1$
correspondence between global minimizers of $I_{q}$ and minimizers of
$\int_{\Omega}|\nabla u|^{2}$ over the set
\[
\left\{  u \in H_{0}^{1}(\Omega): \int_{\Omega}a|u|^{q+1}=1\right\}  .
\]
By \cite[Theorem 1.1]{KLP}, we infer that if $U_{q}$ and $V_{q}$ are global
minimizers of $I_{q}$ then $U_{q}=tV_{q}$ for some $t>0$. But since $U_{q}$
and $V_{q}$ solve $(P_{a,q})$, we
deduce that $t=1$, i.e. $U_{q}$ is the unique nonnegative global minimizer of
$I_{q}$.

\begin{enumerate}
\item Assume by contradiction that $U_{q}(x)=0$ for some $x\in\Omega_{+}$.
Then, by the strong maximum principle, $U_{q}$ vanishes is some ball
$B\subset\Omega_{+}$. We choose a nontrivial $\psi\geq0$ such that $\psi\in
C_{c}^{1}(B)$ and extend it by zero to $\Omega$. Then
\[
I_{q}(U_{q}+t\psi)=I_{q}(U_{q})+I_{q}(t\psi),
\]
and
\[
I_{q}(t\psi)=\frac{t^{2}}{2}\int_{\Omega}|\nabla\psi|^{2}-\frac{t^{q+1}}%
{q+1}\int_{\Omega}a\psi^{q+1}<0
\]
if $t$ is small enough. Hence
\[
I_{q}(U_{q}+t\psi)<I_{q}(U_{q}),
\]
which is a contradiction.

Now, let $q_{0}\in(0,1)$. We will show that $\displaystyle\lim_{q\rightarrow
q_{0}}U_{q}=U_{q_{0}}$ in $W_{D}^{2,r}(\Omega)$. Let $q_{n}\rightarrow q_{0}$.
Then $U_{q_{n}}$ is bounded in $W_{D}^{2,r}(\Omega)$ and we can show that, up
to a subsequence, $U_{q_{n}}\rightarrow\overline{U}$ in $C_{0}^{1}%
(\overline{\Omega})$ for some $\overline{U}$. Since $U_{q_{n}}\rightharpoonup
\overline{U}$ in $H_{0}^{1}(\Omega)$ and $U_{q_{n}}\rightarrow\overline{U}$ in
$L^{p}(\Omega)$ for $1\leq p<2^{\ast}$, we have
\[
I_{q_{0}}(\overline{U})\leq\liminf I_{q_{n}}(U_{q_{n}}).
\]
On the other hand, since $U_{q_{n}}$ is the global minimizer of $I_{q_{n}}$,
we get
\[
\liminf I_{q_{n}}(U_{q_{n}})\leq\liminf I_{q_{n}}(U_{q_{0}})=I_{q_{0}%
}(U_{q_{0}}).
\]
Hence $\overline{U}$ is a ground state solution of $(P_{a,q_{0}})$ so that, by
uniqueness, $\overline{U}=U_{q_{0}}$. Finally, Lemma \ref{lem:c} (i) gives
that $U_{q_{n}}\rightarrow U_{q_{0}}$ in $W_{D}^{2,r}(\Omega)$, as desired.

\item We argue as in the proof of Theorem 1.3 in \cite{krqu}. Let $B$ be a
ball such that $\overline{B}\subset\Omega_{+}$.
Since $(P_{a,q})$ is homogeneous, by a rescaling argument we can assume that
$\lambda_{1}(a,B)<1$. Assume by contradiction that $q_{n}\rightarrow1^{-}$ and
$u_{n}:=U_{q_{n}}\not \in \mathcal{P}^{\circ}$. If $\left\{  u_{n}\right\}  $
is bounded in $H_{0}^{1}(\Omega)$ then we can assume that $u_{n}\rightarrow
u_{0}$ in $H_{0}^{1}(\Omega)$ and $u_{0}$ is a weak solution of
\begin{equation}%
\begin{cases}
-\Delta u_{0}=a\left(  x\right)  u_{0} & \mbox{ in }\Omega,\\
u_{0}=0 & \mbox{ on }\partial\Omega.
\end{cases}
\label{au}%
\end{equation}
Since $u_{n}>0$ in $\Omega_{+}$ we have, by Lemma 2.5 in \cite{krqu},
that there exists some $\phi$ such that $u_{n}\geq\phi>0$ in $B$ for every
$n$, and consequently $u_{0}\not \equiv 0$. The rest of the proof is carried
out as the one of Theorem 1.3 in \cite{krqu}.

\item Recall that
\[
I_{0}(u):=\frac{1}{2}\int_{\Omega}|\nabla u|^{2}-\int_{\Omega}a|u|
\]
for $u\in H_{0}^{1}(\Omega)$. We denote by $U_{0}$ a global minimizer of
$I_{0}$. Let $q_{n}\rightarrow0^{+}$ and $u_{n}:=U_{q_{n}}$. Then $I_{q_{n}%
}(u_{n})\leq I_{q_{n}}(U_{0})$, so that $\liminf I_{q_{n}}(u_{n})\leq
I_{0}(U_{0})$. On the other hand, since $\left\{  u_{n}\right\}  $ is bounded
in $H_{0}^{1}(\Omega)$, we can assume that $u_{n}\rightharpoonup u_{0}$ in
$H_{0}^{1}(\Omega)$, $u_{n}\rightarrow u_{0}$ in $L^{p}(\Omega)$ for $1\leq
p<2^{\ast}$, and $u_{n}\rightarrow u_{0}$ a.e. in $\Omega$. In particular
$u_{0}\geq0$. Hence
\[
I_{0}(u_{0})\leq\liminf I_{q_{n}}(u_{n})\leq I_{0}(U_{0}),
\]
which implies that $u_{0}$ is also a global minimizer of $I_{0}$. Note now
that $\mathcal{S}(a)$ is the unique global minimizer of $\tilde{I}$, given by
\[
\tilde{I}(u):=\frac{1}{2}\int_{\Omega}|\nabla u|^{2}-\int_{\Omega}au,
\]
for $u\in H_{0}^{1}(\Omega)$. Thus, if $\mathcal{S}(a)\geq0$ then
\[
\tilde{I}(u_{0})=I_{0}(u_{0})\leq I_{0}(\mathcal{S}(a))=\tilde{I}%
(\mathcal{S}(a)),
\]
i.e. $u_{0}$ is also a global minimizer of $\tilde{I}$. Therefore $u_{0}%
\equiv\mathcal{S}(a)$. Finally, taking $u_{n}-\mathcal{S}(a)$ as test function
in $(P_{a,q_{n}})$ we get that
\[
\int_{\Omega}\nabla u_{n}\nabla(u_{n}-\mathcal{S}(a))=\int_{\Omega}%
au_{n}^{q_{n}}(u_{n}-\mathcal{S}(a))\rightarrow0,
\]
where we used that $u_{n}^{q_{n}}\leq\max\{1,u_{n}\}$. It follows that
$\int_{\Omega}|\nabla u_{n}|^{2}\rightarrow\int_{\Omega}|\nabla\mathcal{S}%
(a)|^{2}$, and consequently $u_{n}\rightarrow\mathcal{S}(a)$ in $H_{0}%
^{1}(\Omega)$. Using Remark 3.1 in \cite{krqu}, we get $u_{n}\rightarrow
\mathcal{S}(a)$ in $C_{0}^{1}(\overline{\Omega})$,
as desired.
Finally, if $\mathcal{S}(a)>0$ in $\Omega$, then the last
assertion of item (iii) follows from Lemma \ref{lem:c} (ii). $\blacksquare$


\end{enumerate}


\begin{proposition}
\label{pgs2}\strut

\begin{enumerate}
\item If $\lambda_{1}(a)>1$ then $U_{q}\rightarrow0$ in $C_{0}^{1}%
(\overline{\Omega})$ as $q\rightarrow1^{-}$.

\item If $\lambda_{1}(a)<1$ then $\Vert U_{q}\Vert_{\infty}\rightarrow\infty$
as $q\rightarrow1^{-}$.
\end{enumerate}
\end{proposition}

\textit{Proof}. Let $q_{n} \to1^{-}$ and $u_{n}:=U_{q_{n}}$.

\begin{enumerate}
\item First we show that $\left\{  u_{n}\right\}  $ is bounded in $H_{0}%
^{1}(\Omega)$. If not then we can assume that $\Vert u_{n}\Vert\rightarrow
\infty$ and $v_{n}:=\frac{u_{n}}{\Vert u_{n}\Vert}\rightharpoonup v_{0}$ in
$H_{0}^{1}(\Omega)$, $v_{n}\rightarrow v_{0}$ in $L^{p}(\Omega)$ for $1\leq
p<2^{\ast}$ and $v_{n}\rightarrow v_{0}$ a.e. in $\Omega$. Note that
$v_{n}\geq0$ satisfies
\[
-\Delta v_{n}=a(x)\frac{v_{n}^{q_{n}}}{\Vert u_{n}\Vert^{1-q_{n}}}.
\]
Hence
\begin{equation}
1=\int_{\Omega}|\nabla v_{n}|^{2}=\frac{1}{\Vert u_{n}\Vert^{1-q_{n}}}%
\int_{\Omega}a(x)v_{n}^{q_{n}+1}, \label{uno}%
\end{equation}
so that $\Vert u_{n}\Vert^{1-q_{n}}$ is bounded. Since $\Vert u_{n}\Vert\geq1$
for $n$ large enough, we can assume that $\Vert u_{n}\Vert^{1-q_{n}%
}\rightarrow d\geq1$. It follows that $v_{0}\in H_{0}^{1}(\Omega)$ is a
nonnegative solution of
\[
-\Delta v_{0}=\frac{1}{d}a(x)v_{0}.
\]
Finally, $v_{0}\not \equiv 0$, since otherwise $\int_{\Omega}a\left(
x\right)  v_{n}^{q_{n}+1}\rightarrow0$, which contradicts (\ref{uno}). By the
strong maximum principle, we deduce that $v_{0}\in\mathcal{P}^{\circ}$. It
follows that $\lambda_{1}(a)=\frac{1}{d}\leq1$, a contradiction. Hence
$\{u_{n}\}$ is bounded, and arguing as in the proof of Theorem 1.3 in
\cite{krqu}, we can show that, up to a subsequence, $u_{n}\rightarrow u_{0}$
in $C^{1}(\overline{\Omega})$ and $u_{0}\geq0$ solves \eqref{au}. If
$u_{0}\not \equiv 0$ then, by the strong maximum principle, we have $u_{0}%
\in\mathcal{P}^{\circ}$, so that $\lambda_{1}(a)=1$, and we reach a
contradiction again. Therefore $u_{0}\equiv0$. By standard elliptic
regularity, we infer that $u_{n}\rightarrow0$ in $C^{1}(\overline{\Omega})$.

\item It is enough to show that any subsequence of $\{u_{n}\}$ is unbounded in
$H_{0}^{1}(\Omega)$. Assume by contradiction that $\{u_{n}\}$ has a bounded
subsequence in $H_{0}^{1}(\Omega)$, still denoted by $\{u_{n}\}$. By the final
argument in the previous item, up to a subsequence, we have $u_{n}\rightarrow
u_{0}$ in $C^{1}(\overline{\Omega})$, and $u_{0}\geq0$ is a solution of
\eqref{au}. Let us show that $u_{0}\not \equiv 0$, in which case $\lambda
_{1}(a)=1$, and we get a contradiction. Let $\phi>0$ be an eigenfunction
associated to $\lambda_{1}(a)<1$, so that
\[
\int_{\Omega}|\nabla\phi|^{2}-\int_{\Omega}a(x)\phi^{2}<\int_{\Omega}%
|\nabla\phi|^{2}-\lambda_{1}(a)\int_{\Omega}a(x)\phi^{2}=0.
\]
Thus
\begin{align*}
\frac{1}{2}\left(  \int_{\Omega}|\nabla u_{0}|^{2}-\int_{\Omega}a(x)u_{0}%
^{2}\right)   &  =\lim I_{q_{n}}(u_{n})\leq\lim I_{q_{n}}(\phi)\\
&  =\frac{1}{2}\left(  \int_{\Omega}|\nabla\phi|^{2}-\int_{\Omega}a(x)\phi
^{2}\right)  <0,
\end{align*}
which shows that $u_{0}\not \equiv 0$, and the proof is complete.
$\blacksquare$\newline
\end{enumerate}



\section{An implicit function theorem approach I}

\label{sec:b1}

In this section, we discuss the existence of positive solutions for
$(P_{a,q})$ using the implicit function theorem. To this end, we consider the
nonlinear mapping $\mathcal{N}(q,u)=a(x)u^{q}$ for $u\in W_{D}^{2,t}(\Omega)$,
$t>N$, such that $u\geq0$. Note that the Fr\'{e}chet derivative of
$\mathcal{N}$ with respect to $u$ is formally given by
\[
\mathcal{N}_{u}(q,u)h=a(x)qu^{q-1}h,
\]
which is not well-defined in general, since $u=0$ on $\partial\Omega$ and
$q-1<0$. To overcome this difficulty, we shall additionally impose a decay
condition on $a$ near $\partial\Omega$ and take $u$ in $\mathcal{P}^{\circ
}\cap W_{D}^{2,t}(\Omega)$, which is an open subset of $W_{D}^{2,t}(\Omega)$.

Let $(q_{0},u_{0})\in\lbrack0,1)\times W_{D}^{2,t}(\Omega)$, $t>N$, be such
that $u_{0}\in\mathcal{P}^{\circ}$ is a solution of $(P_{a,q_{0}})$. We
consider the nonlinear mapping
\[
\mathcal{F}:U_{0}:=\left(  q_{0}-\frac{\sigma_{0}}{2},q_{0}+\frac{\sigma_{0}%
}{2}\right)  \times B_{0}\rightarrow L^{t}(\Omega);\quad\mathcal{F}%
(q,u):=-\Delta u-a(x)u^{q},
\]
where $\sigma_{0}>0$ and $B_{0}$ is an open ball in $W_{D}^{2,t}(\Omega)$,
centered at $u_{0}$. Let us show how we shall fix $\sigma_{0}$, $B_{0}$ and
$t$. Recall that $\alpha$ and $\rho_{0}$ are given by $(H_{1})$. Since
$W_{D}^{2,t}(\Omega)\subset C_{0}^{1}(\overline{\Omega})$, we may choose
$B_{0}\subset\mathcal{P}^{\circ}$. Furthermore, we pick $B_{0}$,
$0<c_{1}<c_{2}$ and $\rho_{1}\in(0,\min\left\{  1,\rho_{0}\right\}  )$ such
that
\begin{equation}
c_{1}d(x,\partial\Omega)\leq u(x)\leq c_{2}d(x,\partial\Omega)\quad
\mbox{for }x\in\Omega_{\rho_{1}},\text{ and }u\in B_{0}. \label{c1c2}%
\end{equation}
On the other hand, from $(H_{1})$, we take $\sigma_{0}>0$ small enough such
that
\[
\alpha>1-\frac{1}{N}+\sigma_{0}.
\]
Since $q>q_{0}-\frac{\sigma_{0}}{2}\geq-\frac{\sigma_{0}}{2}$, we deduce that
\[
N(\alpha+q-1)>N\left(  -\frac{1}{N}+\sigma_{0}-\frac{\sigma_{0}}{2}\right)
=-1+\frac{\sigma_{0}N}{2}.
\]
This inequality enables us to take $t\in(N,r)$ depending only on $\sigma_{0}$
and $N$, and such that
\begin{equation}
t(\alpha+q-1)>-1+\frac{\sigma_{0}N}{4}. \label{def:t}%
\end{equation}
We have thus fixed $\sigma_{0}$, $B_{0}$ and $t$. We remark that
$d(\cdot,\partial\Omega)^{-1+\frac{\sigma_{0}N}{4}}\in L^{1}(\Omega_{\rho_{1}%
})$.

Under these conditions, $\mathcal{F}$ and its Fr\'{e}chet derivative
$\mathcal{F}_{u}(q,u)$ are well defined. More precisely, $\mathcal{F}$ maps
$U_{0}$ continuously into $L^{t}(\Omega)$. Indeed, from $(H_{1})$,
\eqref{c1c2} and \eqref{def:t}, it follows that
\[
|a\left(  x\right)  u^{q}(x)|^{t}\leq Cd(x,\partial\Omega)^{t(\alpha+q)}\leq
Cd(x,\partial\Omega)^{-1+\frac{\sigma_{0}N}{4}}\quad\text{for a.e. }x\in
\Omega_{\rho_{1}}.
\]
So $au^{q}\in L^{t}(\Omega_{\rho_{1}})$.
Since $0<u\in C(\overline{\Omega})$ and $a\in L^{r}(\Omega)$, it follows that
$au^{q}\in L^{t}(\Omega)$. Hence, $\mathcal{F}$ maps $U_{0}$ into
$L^{t}(\Omega)$. To verify that $\mathcal{F}$ is continuous, let $(q_{n}%
,u_{n})\rightarrow(q_{0},u_{0})$ in $U_{0}$. Then, $u_{n}\rightarrow u_{0}$ in
$C_{0}^{1}(\overline{\Omega})$ and $u_{0}>0$ in $\Omega$, so
that $au_{n}^{q_{n}}\rightarrow au_{0}^{q_{0}}$ a.e. in $\Omega_{\rho_{1}}$.
Moreover, we deduce from \eqref{def:t} that for a.e. $x\in\Omega_{\rho_{1}}$,
\begin{align*}
|a\left(  x\right)  u_{n}^{q_{n}}(x)-a\left(  x\right)  u_{0}^{q_{0}}(x)|^{t}
&  \leq C\left(  d(x,\partial\Omega)^{t(\alpha+q_{n})}+d(x,\partial
\Omega)^{t(\alpha+q_{0})}\right) \\
&  \leq2Cd(x,\partial\Omega)^{-1+\frac{\sigma_{0}N}{4}}\in L^{1}(\Omega
_{\rho_{1}}).
\end{align*}
The Lebesgue dominated convergence theorem shows that $au_{n}^{q_{n}%
}\rightarrow au_{0}^{q_{0}}$ in $L^{t}(\Omega_{\rho_{1}})$. In a similar
manner as above, the desired assertion follows. Next, we formally infer that,
for $(q,u)\in U_{0}$,
\begin{equation}
\mathcal{F}_{u}(q,u)\phi=-\Delta\phi-qa(x)u^{q-1}\phi\label{F_u}%
\end{equation}
(the exact deduction of \eqref{F_u} will be developed in the proof of
Proposition \ref{prop:open} below). Using $(H_{1})$, \eqref{c1c2} and
\eqref{def:t} again, we observe that
\[
|qa\left(  x\right)  u\left(  x\right)  ^{q-1}|^{t}\leq Cd(x,\partial
\Omega)^{t(\alpha+q-1)}\leq Cd(x,\partial\Omega)^{-1+\frac{\sigma_{0}N}{4}%
}\quad\text{for a.e.\ }x\in\Omega_{\rho_{1}}.
\]
We deduce then (in the same way as for $\mathcal{F}$) that $\mathcal{F}%
_{u}(q,u)$ is a bounded linear operator from $W_{D}^{2,t}(\Omega)$ to
$L^{t}(\Omega)$.

\begin{proposition}
\label{prop:open} Assume $(H_{1})$, and let $u_{0}\in\mathcal{P}^{\circ}$ be a
solution of $(P_{a,q_{0}})$ with $q_{0}\in\lbrack0,1)$. Then, the Fr\'{e}chet
derivative $\mathcal{F}_{u}(q_{0},u_{0})$ at $(q_{0},u_{0})$ maps $W_{D}%
^{2,t}(\Omega)$ onto $L^{t}(\Omega)$ \textrm{homeomorphically}, and there
exists a curve $q\mapsto u(q)$ from $(q_{0}-\delta_{0},q_{0}+\delta_{0})$ to
$W_{D}^{2,t}(\Omega)$, for some $\delta_{0}>0$, such that $u(q_{0})=u_{0}$,
$F(q,u(q))=0$, and $u(q)\in\mathcal{P}^{\circ}$ for $(q_{0}-\delta_{0}%
,q_{0}+\delta_{0})$. In particular, $\mathcal{I}_{a}$ is open.
\end{proposition}

\textit{Proof.} We verify that $\mathcal{F}_{u}(q,u)$ is well defined for
$(q,u)\in U_{0}$ as a bounded linear operator from $W_{D}^{2,t}(\Omega)$ to
$L^{t}(\Omega)$. Set $\mathcal{N}(q,u):=a(x)u^{q}$, and consider the
F\'{r}echet derivative $\mathcal{N}_{u}(q,u)$. Using the mean value theorem,
we find $\theta\in(0,1)$ such that
\[
\mathcal{N}(q,u+h)-\mathcal{N}(q,u)=qa(x)(u+\theta h)^{q-1}h,
\]
where $\Vert h\Vert_{2,t}$ is small enough so that $u+h,u+\theta h\in B_{0}$.
It follows that
\[
\frac{\Vert\mathcal{N}(q,u+h)-\mathcal{N}(q,u)-qa(x)u^{q-1}h\Vert_{t}}{\Vert
h\Vert_{2,t}}=\frac{\Vert qa(x)h\{(u+\theta h)^{q-1}-u^{q-1}\}\Vert_{t}}{\Vert
h\Vert_{2,t}}.
\]
We know that $h\in C(\overline{\Omega})$ and $\Vert h\Vert_{C(\overline
{\Omega})}\leq C\Vert h\Vert_{2,t}$. Using these facts, we deduce that
\[
\frac{\Vert qa(x)h\{(u+\theta h)^{q-1}-u^{q-1}\}\Vert_{t}}{\Vert h\Vert_{2,t}%
}\leq C\Vert a(x)\{(u+\theta h)^{q-1}-u^{q-1}\}\Vert_{t}.
\]
Now we use the Lebesgue dominated convergence theorem to show that
\begin{equation}
\Vert a(x)\{(u+\theta h)^{q-1}-u^{q-1}\}\Vert_{t}\rightarrow0\quad
\mbox{ as }\ \Vert h\Vert_{2,t}\rightarrow0. \label{LCT}%
\end{equation}
Indeed, since $u+\theta h\in B_{0}$,
\eqref{c1c2} and \eqref{def:t} imply that, for a.e.\ $x\in\Omega_{\rho_{1}}$,
\begin{align*}
|a(x)\{(u+\theta h)^{q-1}-u^{q-1}\}|^{t}  &  \leq C\{|a(x)(u+\theta
h)^{q-1}|^{t}+|a(x)u^{q-1}|^{t}\}\\
&  \leq C^{\prime}d(x,\partial\Omega)^{t(\alpha+q-1)}\\
&  \leq C^{\prime}d(x,\partial\Omega)^{-1+\frac{\sigma_{0}N}{4}}\in
L^{1}(\Omega_{\rho_{1}}).
\end{align*}
For a.e. $x\in\Omega\setminus\Omega_{\rho_{1}}$, we have that
\[
|a(x)\{(u+\theta h)^{q-1}-u^{q-1}\}|^{t}\leq C|a(x)|^{t}\in L^{1}(\Omega).
\]
Moreover, if $\Vert h\Vert_{2,t}\rightarrow0$, then
\[
a(x)\{(u+\theta h)^{q-1}-u^{q-1}\}\rightarrow0\quad\mbox{a.e. in }\Omega.
\]
The Lebesgue dominated convergence theorem can thus be applied to deduce \eqref{LCT}.

Therefore, we have obtained that
\[
\lim_{\Vert h\Vert_{2,t}\rightarrow0}\frac{\Vert\mathcal{N}(q,u+h)-\mathcal{N}%
(q,u)-qa(x)u^{q-1}h\Vert_{t}}{\Vert h\Vert_{2,t}}=0.
\]
In addition, using \eqref{c1c2} and \eqref{def:t} again, we find that the
mapping $h\mapsto qa(x)u^{q-1}h$, from $W_{D}^{2,t}(\Omega)$ to $L^{t}%
(\Omega)$, is linear and bounded. Summing up, we have verified that
\[
\mathcal{N}_{u}(q,u)h=qa(x)u^{q-1}h,\quad\text{and}\quad\mathcal{F}%
_{u}(q,u)h=-\Delta h-qa(x)u^{q-1}h.
\]

Next, we shall show how to apply the implicit function theorem \cite[Theorem
4.B]{Zei86} to $(q_{0},u_{0})$ such that $\mathcal{F}(q_{0},u_{0})=0$, with
$q_{0}\in\lbrack0,1)$ and $u_{0}\in\mathcal{P}^{\circ}$. In a similar manner,
relying on $(H_{1})$, \eqref{c1c2} and \eqref{def:t}, we can check that
$\mathcal{F}_{u}(\cdot,\cdot):U_{0}\rightarrow\mathcal{L}(W_{D}^{2,t}%
(\Omega),L^{t}(\Omega))$ is continuous. We consider first the case $q_{0}>0$.
Note that
\[
\mathcal{F}_{u}(q_{0},u_{0})\phi=-\Delta\phi-q_{0}a(x)u_{0}^{q_{0}-1}\phi.
\]
We claim that
\begin{equation}
\mathcal{F}_{u}(q_{0},u_{0}):W_{D}^{2,t}(\Omega)\rightarrow L^{t}%
(\Omega)\ \ \mbox{is homeomorphic.} \label{homeo}%
\end{equation}
To verify it, we study the eigenvalue problem
\[
\mathcal{F}_{u}(q_{0},u_{0})\phi=\sigma\phi.
\]
By $\sigma_{1}=\sigma_{1}(q_{0},u_{0})$ we denote the smallest eigenvalue
(which is simple) of this equation, and by $\phi_{1}$ a positive eigenfunction
belonging to $\mathcal{P}^{\circ}$, associated to $\sigma_{1}$. Using the
divergence theorem (as stated e.g. in \cite[p. 742]{cuesta}), we can deduce
that
\begin{equation}
\int_{\Omega}(-\Delta u_{0})q_{0}u_{0}^{q_{0}-1}\phi_{1}+u_{0}^{q_{0}}%
\Delta\phi_{1}=\int_{\Omega}|\nabla u_{0}|^{2}q_{0}(q_{0}-1)u_{0}^{q_{0}%
-2}\phi_{1}. \label{div}%
\end{equation}
Indeed, we first note that both sides in \eqref{div} are well defined, since
$u_{0}$ and $\phi_{1}$ behave like $d(x,\partial\Omega)$ in $\Omega_{\rho_{1}%
}$ and $u_{0}$, $\phi_{1}$ are positive in $\Omega$. For instance, for
$|\nabla u_{0}|^{2}$ and $u_{0}^{q_{0}-2}\phi_{1}$, we see that, for
$x\in\Omega_{\rho_{1}}$
\begin{align*}
&  |\nabla u_{0}|^{2}\leq C,\\
&  |u_{0}^{q_{0}-2}\phi_{1}|\leq Cd(x,\partial\Omega)^{q_{0}-1}\in
L^{1}(\Omega_{\rho_{1}}),
\end{align*}
where we have used the fact that $q_{0}-1>-1$. Let us check equality
\eqref{div}. A direct computation yields
\begin{align*}
\int_{\Omega}\mathrm{div}\left(  (\nabla u_{0})q_{0}u_{0}^{q_{0}-1}\phi
_{1}\right)   &  =\int_{\Omega}(\Delta u_{0})q_{0}u_{0}^{q_{0}-1}\phi_{1}\\
&  +\int_{\Omega}q_{0}(q_{0}-1)|\nabla u_{0}|^{2}u_{0}^{q_{0}-2}\phi_{1}\\
&  +\int_{\Omega}q_{0}(\nabla u_{0}\nabla\phi_{1})u_{0}^{q_{0}-1}.
\end{align*}
Since $u_{0}^{q_{0}-1}\phi_{1}\nabla u_{0}\in\left(  W^{1,\gamma_{0}}%
(\Omega)\right)  ^{n}$ for some $\gamma_{0}=\gamma_{0}(q_{0})>1$, the
divergence theorem applies, and we obtain that
\begin{align*}
\int_{\Omega}\mathrm{div}\left(  (\nabla u_{0})q_{0}u_{0}^{q_{0}-1}\phi
_{1}\right)   &  =\int_{\Omega}\mathrm{div}\left(  (\nabla u_{0})q_{0}\left(
\frac{\phi_{1}}{u_{0}}\right)  ^{1-q_{0}}\phi_{1}^{q_{0}}\right) \\
&  =\int_{\partial\Omega}\frac{\partial u_{0}}{\partial\nu}q_{0}\left(
\frac{\phi_{1}}{u_{0}}\right)  ^{1-q_{0}}\phi_{1}^{q_{0}}\\
&  =0,
\end{align*}
where we have used that $\phi_{1}=0$ on $\partial\Omega$. In a similar manner,
we deduce that, by a direct computation,
\[
\int_{\Omega}\mathrm{div}\left(  (\nabla\phi_{1})u_{0}^{q_{0}}\right)
=\int_{\Omega}(\Delta\phi_{1})u_{0}^{q_{0}}+\int_{\Omega}(\nabla\phi_{1}\nabla
u_{0})q_{0}u_{0}^{q_{0}-1},
\]
and, by the divergence theorem,
\[
\int_{\Omega}\mathrm{div}\left(  (\nabla\phi_{1})u_{0}^{q_{0}}\right)
=\int_{\partial\Omega}\frac{\partial\phi_{1}}{\partial\nu}u_{0}^{q_{0}}=0,
\]
where we have used that $u_{0}=0$ on $\partial\Omega$. Combining these
assertions we obtain \eqref{div}.

From \eqref{div}, it follows that
\begin{align*}
\int_{\Omega}|\nabla u_{0}|^{2} q_{0}(q_{0} - 1) u_{0}^{q_{0}-2} \phi_{1}  &
= \int_{\Omega}(-\Delta u_{0}) q_{0} u_{0}^{q_{0} -1} \phi_{1} + u_{0}^{q_{0}}
(\Delta\phi_{1})\\
&  = \int_{\Omega}(a u_{0}^{q_{0}}) q_{0} u_{0}^{q_{0} -1} \phi_{1} +
u_{0}^{q_{0}} (-q_{0}au_{0}^{q_{0}-1}\phi_{1} - \sigma_{1} \phi_{1})\\
&  = -\sigma_{1} \int_{\Omega}u_{0}^{q_{0}} \phi_{1},
\end{align*}
and thus that
\begin{align*}
\sigma_{1} = \frac{q_{0}(1- q_{0}) \int_{\Omega}|\nabla u_{0}|^{2}
u_{0}^{q_{0}-2} \phi_{1} }{\int_{\Omega}u_{0}^{q_{0}} \phi_{1}} > 0.
\end{align*}
The assertion $\sigma_{1} > 0$ tells us that $\mathcal{F}_{u}(q_{0}, u_{0})$
is bijective. Since $\mathcal{F}_{u}(q_{0}, u_{0})$ is continuous, the Bounded
Inverse Theorem yields \eqref{homeo}.

It
remains to consider the case $q_{0}=0$. However, we note that $\mathcal{F}%
_{u}(0, u_{0})=-\Delta$ maps $W^{2,t}_{D}(\Omega)$ onto $L^{t}(\Omega)$ homeomorphically.

We are now ready to apply \cite[Theorem 4.B]{Zei86} to $\mathcal{F}$ at
$(q_{0},u_{0})$, which provides us with some $\delta_{0}>0$ such that
$\mathcal{F}(q,u(q))=0$ for $q\in(q_{0}-\delta_{0},q_{0}+\delta_{0})$,
$q\mapsto u(q)\in W_{D}^{2,t}(\Omega)$ is continuous, and $u(q_{0})=u_{0}$. In
particular, $q\mapsto u(q)\in C^{1}(\overline{\Omega})$ is continuous, so that
$u(q)\in\mathcal{P}^{\circ}$ for every $q\in(q_{0}-\delta_{0},q_{0}+\delta
_{0})$, since $u(q_{0})\in\mathcal{P}^{\circ}$. In particular, $(q_{0}%
-\delta_{0},q_{0}+\delta_{0})\subset\mathcal{I}_{a}$, as desired.
$\blacksquare$

\begin{remark}
\label{si} \strut

\begin{enumerate}
\item By the uniqueness of positive solutions, the set of positive solutions
of $(P_{a,q})$, with $q$ close to $1$, consists of a curve $\{(q,u(q))\in
\left(  q_{1},1\right)  \times\mathcal{P}^{\circ}\}$ in $\mathbb{R}\times
W_{D}^{2,t}(\Omega)$, for some $0<q_{1}<1$.

\item Using the implicit function theorem \cite[Theorem 4.B(d)]{Zei86}, we can
deduce that $q\mapsto u(q)$ is $C^{1}$ around $q_{0}\in\lbrack0,1)$. Formally,
we show that
\begin{equation}
\mathcal{F}_{q}(q,u)=-a\left(  x\right)  u^{q}\log u. \label{F_q}%
\end{equation}
Indeed, given $\sigma\in(0,1)$ and $s_{0}>0$, we have $|\log s|\leq
Cs^{-\sigma}$ for $0<s\leq s_{0}$ and some $C>0$. Hence, we deduce from
$(H_{1})$, \eqref{c1c2} and \eqref{def:t} that
\[
|a\left(  x\right)  u^{q}\log u|^{t}\leq Cd(x,\partial\Omega)^{t(\alpha
+q-\sigma)}\leq Cd(x,\partial\Omega)^{-1+\frac{\sigma_{0}N}{4}}\quad
\mbox{for }\text{a.e. }x\in\Omega_{\rho_{1}}.
\]
This implies that $a\left(  x\right)  u^{q}\log u\in L^{t}(\Omega_{\rho_{1}}%
)$. By the same argument as for $\mathcal{F}_{u}(q,u)$, we obtain $a\left(
x\right)  u^{q}\log u\in L^{t}(\Omega)$, so that \eqref{F_q} is proved,
$\mathcal{F}_{q}(\cdot,\cdot)$ is continuous (by use of the Lebesgue dominated
convergence theorem), and thus $\mathcal{F}$ is $C^{1}$. The conclusion now
follows.\newline
\end{enumerate}
\end{remark}

By Lemma \ref{lem:c}, we have the following stronger result on the continuity
of $u(q)$ with respect to $q\in(q_{0}-\delta_{0},q_{0}+\delta_{0})$:

\begin{corollary}
\label{cor:r} Under the assumptions of Proposition \ref{prop:open}, we have
the following:

\begin{enumerate}
\item If $q_{0}\in(0,1)$, then the mapping $q\mapsto u(q)$ is continuous from
$(q_{0}-\delta_{0},q_{0}+\delta_{0})$ to $W_{D}^{2,r}(\Omega)$.

\item If $q_{0}=0$, then it is continuous from $[0,\delta_{0})$ to
$W_{D}^{2,r}(\Omega)$.
\end{enumerate}
\end{corollary}

\textit{Proof.} We prove assertions (i) and (ii), using Lemma \ref{lem:c} (i)
and (ii), respectively. For assertion (i), let $q_{1}\in(q_{0}-\delta
_{0},q_{0}+\delta_{0})$. Then, we know that $u(q)\rightarrow u(q_{1})$ in
$W_{D}^{2,t}(\Omega)$ as $q\rightarrow q_{1}$, and $W_{D}^{2,t}(\Omega)\subset
C(\overline{\Omega})$, so that $u(q)\rightarrow u(q_{1})$ in $C(\overline
{\Omega})$. Assertion (i) is now a direct consequence of Lemma \ref{lem:c}
(i). The proof of item (ii) is similar, so we omit it. $\blacksquare$ \newline


\section{An implicit function theorem approach II}

\label{sec:b2}

Excepting Corollary \ref{cor:0to1}, throughout this section we assume that
$\lambda_{1}(a)=1$, so that $(P_{a,q})$ possesses the trivial line of
solutions%
\[
\Gamma_{1}:=\left\{  (q,u)=(1,s\phi_{1}):s>0\right\}  ,
\]
where $\phi_{1}$ is the positive eigenfunction ($\Vert\phi_{1}\Vert_{2}=1$)
associated with $\lambda_{1}(a)=1$.
We shall look at $q$ as a bifurcation parameter in $(P_{a,q})$, and then
seek for bifurcating solutions in $\mathcal{P}^{\circ}$ from the trivial line
$\Gamma_{1}$. To this end, we employ the Lyapunov-Schmidt reduction, based on
the positive eigenfunction $\phi_{1}$. We will construct solutions of
$(P_{a,q})$ bifurcating from a certain point $(1,t\phi_{1})\in\Gamma_{1}$, in
the topology $\mathbb{R}\times W_{D}^{2,\eta}(\Omega)$ for some $\eta>N$.
Consequently they also belong to $\mathcal{P}^{\circ}$.

Now,
we weaken the decay condition on $a$ used in Section \ref{sec:b1}. Namely, we
assume:
\[
\left\vert a(x)\right\vert \leq Cd(x,\partial\Omega)^{\alpha}\text{
\ \negthinspace a.e.\ in }\Omega_{\rho_{0}},\text{ for some }\rho_{0}>0\text{
and }\alpha>-\frac{1}{N}.\leqno{(H_1')}
\]
We pick $\sigma_{0}>0$ small enough such that $\alpha>\sigma_{0}-\frac{1}{N}$,
and set $I_{0}:=(1-\frac{\sigma_{0}}{2},1+\frac{\sigma_{0}}{2})$. We see that
there exists $\eta\in(N,r)$ such that
\begin{equation}
\eta(\alpha+q-1)>-1+\frac{\sigma_{0}N}{4} \label{Lrcond}%
\end{equation}
for $q\in I_{0}$. Note that $\eta$ can be determined depending only on $N$ and
$\sigma_{0}$. In the sequel, we fix $\eta$ in this way.

We set
\[
A:=-\Delta-a(x)\quad\mbox{and}\quad D(A):=W_{D}^{2,\eta}(\Omega).
\]
It follows that $\mathrm{Ker}A=\langle\phi_{1}\rangle:=\{s\phi_{1}%
:s\in\mathbb{R}\}$. We split $D(A)$ as follows:
\[
D(A)=\mathrm{Ker}A+X_{2};\quad u=t\phi_{1}+w,
\]
where $t:=\int_{\Omega}u\phi_{1}$, and $w:=u-(\int_{\Omega}u\phi_{1})\phi_{1}%
$. So, $X_{2}$ is characterized as
\[
X_{2}=\left\{  w\in W_{D}^{2,\eta}(\Omega):\int_{\Omega}w\phi_{1}=0\right\}
.
\]
On the other hand, put $Y:=L^{\eta}(\Omega)=Y_{1}+R(A)$, where
\[
Y_{1}=\mathrm{Ker}A,\quad\mbox{ and }\ R(A)=\left\{  f\in L^{\eta}%
(\Omega):\int_{\Omega}f\phi_{1}=0\right\}  .
\]
Let $Q$ be the projection of $Y$ to $R(A)$, given by
\[
Q[f]:=f-\left(  \int_{\Omega}f\phi_{1}\right)  \phi_{1}.
\]
We thus reduce $(P_{a,q})$ to the following coupled equations:
\begin{align*}
&  Q[Au]=Q[a\left(  x\right)  \left(  u^{q}-u\right)  ],\\
&  (1-Q)[Au]=(1-Q)[a\left(  x\right)  \left(  u^{q}-u\right)  ].
\end{align*}
The first equation yields
\begin{equation}
Aw=Q[a\left(  x\right)  \{(t\phi_{1}+w)^{q}-(t\phi_{1}+w)\}], \label{beq01}%
\end{equation}
where we have used the fact that $\int_{\Omega}(\phi_{1}Au-uA\phi_{1})=0$ and
$A\phi_{1}=0$ (and so, $Au=Aw$). The second equation implies that
\begin{align*}
0  &  =(1-Q)[a\left(  x\right)  (u^{q}-u)]\\
&  =\left(  \int_{\Omega}a(x)\left\{  (t\phi_{1}+w)^{q}-(t\phi_{1}+w)\right\}
\phi_{1}\right)  \phi_{1},
\end{align*}
and thus, that
\begin{equation}
0=\int_{\Omega}a(x)\left\{  (t\phi_{1}+w)^{q}-(t\phi_{1}+w)\right\}  \phi_{1}.
\label{beq02}%
\end{equation}

Now, we see that $(q,t,w)=(1,t,0)$ satisfies (\ref{beq01}) and (\ref{beq02})
for any $t>0$. So, first we solve \eqref{beq01} with respect to $w$, around
$(q,t,w)=(1,t_{0},0)$ for a fixed $t_{0}>0$. To this end, we introduce the
mapping
\[
F:I_{0}\times(t_{0}-d,t_{0}+d)\times B_{\rho}(0)\rightarrow R(A)
\]
given by
\[
F(q,t,w):=Aw-Q[a\left(  x\right)  \{(t\phi_{1}+w)^{q}-(t\phi_{1}+w)\}],
\]
where $B_{\rho}(w)$ is the ball in $X_{2}$ centered at $w$, with radius
$\rho>0$. It is clear that $F(1,t_{0},0)=0$, and the condition $\phi_{1}%
\in\mathcal{P}^{\circ}$ tells us that $t\phi_{1}+w\in\mathcal{P}^{\circ}$.
Also, there exist $0<c_{1}<c_{2}$ such that
\[
c_{1}d(x,\partial\Omega)\leq t\phi_{1}+w\leq c_{2}d(x,\partial\Omega
),\quad\text{for }x\in\Omega_{\delta},
\]
where $\sigma_{0},d,\rho$ and $\delta$ are chosen smaller if necessary.
Therefore, thanks to \eqref{Lrcond}, the Fr\'{e}chet derivative $F_{w}%
(q,t,w):X_{2}\rightarrow R(A)$ can be defined by
\begin{equation}
F_{w}(q,t,w)\varphi=A\varphi-Q[a\left(  x\right)  (q(t\phi_{1}+w)^{q-1}%
-1)\varphi], \label{Fw}%
\end{equation}
and moreover, $F_{w}$ is continuous around $(1,t_{0},0)$. We see that
\[
F_{w}(1,t_{0},0)\varphi=A\varphi,
\]
so that
\[
F_{w}(1,t_{0},0)\varphi=0\quad\Longleftrightarrow\quad\varphi=c\phi
_{1}\ \mbox{ for some }c>0.
\]
Since $\varphi\in X_{2}$, it follows that $\int_{\Omega}c\phi_{1}^{2}=0$, and
thus $c=0$. This means that $F_{w}(1,t_{0},0)$ is injective. It is also
surjective, for $\int_{\Omega}f\phi_{1}=0$ if and only if there exists
$\varphi$ such that%
\[%
\begin{cases}
A\varphi=f & \mbox{ in }\Omega,\\
\varphi=0 & \mbox{ on }\partial\Omega.
\end{cases}
\]
Since $F_{w}(1,t_{0},0)$ is continuous, from the Bounded Inverse Theorem we
infer that $F_{w}(1,t_{0},0)$ is an isomorphism. Hence, the implicit function
theorem applies, and consequently, we have
\begin{align*}
&  F(q,t,w)=0,\ \text{for }(q,t,w)\simeq(1,t_{0},0)\\
&  \Longleftrightarrow w=w(q,t),\ \text{for }(q,t)\simeq(1,t_{0}%
)\mbox{ and }w(1,t_{0})=0.
\end{align*}
We plug $w(q,t)$ into \eqref{beq02} to get the following bifurcation equation
in $\mathbb{R}^{2}$:
\[
\Phi(q,t):=\int_{\Omega}a(x)\{(t\phi_{1}+w(q,t))^{q}-(t\phi_{1}+w(q,t)\}\phi
_{1}=0,\quad(q,t)\simeq(1,t_{0}).
\]

For our procedure, we check the following properties of $w$ and $\Phi$:

\begin{lemma}
\strut

\begin{enumerate}
\item $w$ is $C^{1}$ around $(1,t_{0})$ for $t_{0}>0$. Moreover,
$w_{q}(q,t)(\cdot),w_{t}(q,t)(\cdot)\in W_{D}^{2,\eta}(\Omega)$ satisfy
\begin{align}
&  Aw_{q}=Q\left[  a(x)\left\{  (t\phi_{1}+w)^{q}\left(  \log(t\phi
_{1}+w)+\frac{qw_{q}}{t\phi_{1}+w}\right)  -w_{q}\right\}  \right]
,\label{Awq}\\
&  Aw_{t}=Q\left[  a(x)\left\{  q(t\phi_{1}+w)^{q-1}(\phi_{1}+w_{t})-(\phi
_{1}+w_{t})\right\}  \right]  . \label{Awt}%
\end{align}

\item $\Phi$ is $C^{1}$ around $(1,t_{0})$ for $t_{0}>0$. Moreover,
\begin{align}
&  \Phi_{q}(q,t)=\int_{\Omega}a\left(  x\right)  \left[  (t\phi_{1}%
+w)^{q}\left\{  \log(t\phi_{1}+w)+\frac{qw_{q}}{t\phi_{1}+w}\right\}
-w_{q}\right]  \phi_{1}.\label{Phiq}\\
&  \Phi_{t}(q,t)=\int_{\Omega}a(x)\left\{  q(t\phi_{1}+w)^{q-1}(\phi_{1}%
+w_{t})-(\phi_{1}+w_{t})\right\}  \phi_{1}. \label{Phit}%
\end{align}

\end{enumerate}
\end{lemma}

\textit{Proof}. First we prove that $w$ is $C^{1}$. To this end, we verify
that $F$ is $C^{1}$. Thanks to \eqref{Lrcond}, we deduce that
\begin{align}
&  F_{t}(q,t,w)=-Q\left[  a(x)\left\{  q(t\phi_{1}+w)^{q-1}\phi_{1}-\phi
_{1}\right\}  \right]  ,\label{Ft}\\
&  F_{q}(q,t,w)=-Q\left[  a(x)(t\phi_{1}+w)^{q}\log(t\phi_{1}+w)\right]  ,
\label{Fq}%
\end{align}
and moreover, $F_{w}$, $F_{t}$ and $F_{q}$ are continuous around $(1,t_{0}%
,0)$, as desired. Hence, the implicit function theorem yields that $w\in
C^{1}$, and so employing (\ref{beq01}) we see that \eqref{Awq} and \eqref{Awt}
hold. Using \eqref{Lrcond} again, we deduce that $\Phi$ is $C^{1}$, and
\eqref{Phiq} and \eqref{Phit} hold. $\blacksquare$



\begin{proposition}
\label{prop:t*} \label{t1} Suppose that $(H_{1}^{\prime})$ holds, and
$(q_{n},u_{n})\in(0,1)\times\mathcal{P}^{\circ}$ are solutions of
$(P_{a,q_{n}})$ such that $(q_{n},u_{n})\rightarrow(1,t\phi_{1})\in\Gamma_{1}$
in $\mathbb{R}\times W_{D}^{2,\eta}(\Omega)$ for some $t>0$. Then $t=t_{\ast}%
$, where
\[
t_{\ast}:=\exp\left[  -\frac{\int_{\Omega}a\left(  x\right)  \phi_{1}^{2}%
\log\phi_{1}}{\int_{\Omega}a\left(  x\right)  \phi_{1}^{2}}\right]  .
\]

\end{proposition}

\textit{Proof}. Since $(q_{n},u_{n})\rightarrow(1,t\phi_{1})$ in
$\mathbb{R}\times W^{2,\eta}(\Omega)$ for some $t>0$, we have $\Phi
_{q}(1,t)=0$ by the implicit function theorem. Since $w(1,t)=0$, it follows
that
\begin{align*}
\Phi_{q}(1,t)  &  =\int_{\Omega}a\left(  x\right)  \left[  (t\phi_{1})\left\{
\log(t\phi_{1})+\frac{w_{q}(1,t)}{t\phi_{1}}\right\}  -w_{q}(1,t)\right]
\phi_{1}\\
&  =t\int_{\Omega}a\left(  x\right)  \phi_{1}^{2}\log(t\phi_{1})\\
&  =t\left\{  (\log t)\int_{\Omega}a\left(  x\right)  \phi_{1}^{2}%
+\int_{\Omega}a\left(  x\right)  \phi_{1}^{2}\log\phi_{1}\right\}  .
\end{align*}
So, the desired conclusion follows. $\blacksquare$
\smallskip

Next, we consider the existence of bifurcating solutions of $(P_{a,q})$ at
$(1,t_{\ast}\phi_{1})$. To this end, we need to change the choice of $\eta$.
We shall assume $(H_{1})$ and choose $\sigma_{0}>0$ such that $\alpha
>1+\sigma_{0}-\frac{1}{N}$, and set $I_{0}:=(1-\frac{\sigma_{0}}{2}%
,1+\frac{\sigma_{0}}{2})$. Thanks to $(H_{1})$, we fix $\eta\in(N,r)$
depending only on $N$ and $\sigma_{0}$, and such that
\[
\eta(\alpha+q-2)>-1+\frac{\sigma_{0}N}{4},
\]
for $q\in I_{0}$.

Recalling $(H_{1})$, we deduce from \eqref{Fw}, \eqref{Ft} and \eqref{Fq} that
$F$ is $C^{2}$.
The implicit function theorem ensures that so is $w$. Then, we obtain from
\eqref{Awq} and \eqref{Awt} that
\begin{align*}
&  Aw_{qq}=Q\left[  a(x)\left\{  (t\phi_{1}+w)^{q}\left(  \left(  \log
(t\phi_{1}+w)+\frac{qw_{q}}{t\phi_{1}+w}\right)  ^{2}\right.  \right.  \right.
\\
&  \left.  \left.  \left.  \qquad\qquad+\frac{2w_{q}+qw_{qq}}{t\phi_{1}%
+w}-\frac{q(w_{q})^{2}}{(t\phi_{1}+w)^{2}}\right)  -w_{qq}\right\}  \right]
,\\
&  Aw_{qt}=Q\left[  a(x)\left\{  q(t\phi_{1}+w)^{q-1}(\phi_{1}+w_{t})\left(
\log(t\phi_{1}+w)+\frac{qw_{q}}{t\phi_{1}+w}\right)  \right.  \right. \\
&  \left.  \left.  \qquad\qquad+(t\phi_{1}+w)^{q}\left(  \frac{\phi_{1}%
+w_{t}+qw_{qt}}{t\phi_{1}+w}-\frac{qw_{q}(\phi_{1}+w_{t})}{(t\phi_{1}+w)^{2}%
}\right)  -w_{qt}\right\}  \right]  ,\\
&  Aw_{tt}=Q\left[  a(x)\left\{  q(q-1)(t\phi_{1}+w)^{q-2}(\phi_{1}+w_{t}%
)^{2}\right.  \right. \\
&  \left.  \left.  \qquad\qquad+q(t\phi_{1}+w)^{q-1}w_{tt}-w_{tt}\right\}
\right]  .
\end{align*}

In this situation, we can obtain the second derivatives of $\Phi$ by virtue of
$(H_{1})$:

\begin{lemma}
\label{lem:Phi2} $\Phi$ is $C^{2}$ around $(1,t_{0})$ for $t_{0}>0$.
Moreover,
\begin{align}
&  \Phi_{qq} = \int_{\Omega}a(x) \left[  (t\phi_{1} + w)^{q} \left\{  \left(
\log(t\phi_{1} + w) + \frac{q w_{q}}{t \phi_{1} + w} \right)  ^{2} \right.
\right. \nonumber\\
&  \left.  \left.  \qquad+ \frac{2w_{q} + q w_{qq}}{t \phi_{1} + w} - \frac{q
(w_{q})^{2}}{(t \phi_{1} + w)^{2}} \right\}  - w_{qq} \right]  \phi
_{1},\label{Phiqq}\\
&  \Phi_{qt} =\int_{\Omega}a\left(  x\right)  \left[  q(t\phi_{1}%
+w)^{q-1}(\phi_{1}+w_{t})\left\{  \log(t\phi_{1}+w)+\frac{qw_{q}}{t\phi_{1}%
+w}\right\}  \right. \nonumber\\
&  \qquad\left.  +(t\phi_{1}+w)^{q}\left\{  \frac{\phi_{1}+w_{t}}{t\phi_{1}%
+w}+\frac{qw_{qt}(t\phi_{1}+w)-qw_{q}(\phi_{1}+w_{t})}{(t\phi_{1}+w)^{2}%
}\right\}  -w_{qt}\right]  \phi_{1}.\label{Phiqt}\\
&  \Phi_{tt} = \int_{\Omega}a(x) \left[  q(q-1)(t \phi_{1} + w)^{q-2}(\phi_{1}
+ w_{t})^{2} + \left\{  q (t \phi_{1} + w)^{q-1} - 1 \right\}  w_{tt} \right]
\phi_{1}. \label{Phitt}%
\end{align}

\end{lemma}

Based on Lemma \ref{lem:Phi2}, the following existence result is proved.

\begin{proposition}
\label{prop:gam2} \label{t2} Suppose that $(H_{1})$ holds. Then, the set of
solutions of $(P_{a,q})$ consists of $\Gamma_{1}\cup\Gamma_{2}$ in a
neighborhood of $(q,u)=(1,t_{\ast}\phi_{1})$ in $\mathbb{R}\times
W_{D}^{2,\eta}(\Omega)$, where
\[
\Gamma_{2}:=\{(q(s),\,t(s)\phi_{1}+w(q(s),t(s))):|s|<s_{0}\},
\]
for some $s_{0}>0$. Here $s \mapsto q(s),t(s)$ are continuous in
$(-s_{0},s_{0})$ and satisfy:

\begin{enumerate}
\item $q(0)=1$ and $t(0)=t_{\ast}$;

\item $q(s)<1$ for $-s_{0}<s<0$, whereas $q(s)>1$ for $0<s<s_{0}$;

\item $t(s)\phi_{1}+w(q(s),t(s))\in\mathcal{P}^{\circ}$ for $|s|<s_{0}$.
\end{enumerate}
\end{proposition}

\textit{Proof}. We use the Morse Lemma \cite[Theorem 4.3.19]{DM13}
to prove the existence of $\Gamma_{2}$. We know from Proposition \ref{t1}
that
\begin{equation}
\Phi_{q}(1,t_{\ast})=0. \label{Phq1t*}%
\end{equation}
By direct observations of \eqref{Phit} and \eqref{Phitt}, we find that
\begin{equation}
\Phi_{tt}(1,t_{\ast})=\Phi_{t}(1,t_{\ast})=0. \label{Pht:tt1t*}%
\end{equation}
Now, we verify that $\Phi_{qt}(1,t_{\ast})>0$. We derive from \eqref{Phiqt}
that
\begin{equation}
\Phi_{qt}(1,t)=\int_{\Omega}a\left(  x\right)  [(\phi_{1}+w_{t}(1,t))\log
(t\phi_{1})+(\phi_{1}+w_{t}(1,t))]\phi_{1}. \label{Pqt}%
\end{equation}
Since $w_{t}(1,t)=0$ from \eqref{Awt}, it follows from \eqref{Pqt} that
\begin{align*}
\Phi_{qt}(1,t)  &  =\int_{\Omega}a\left(  x\right)  [\phi_{1}\log(t\phi
_{1})+\phi_{1}]\phi_{1}\\
&  =\int_{\Omega}a\left(  x\right)  \phi^{2}\log(t\phi_{1})+\int_{\Omega
}a\left(  x\right)  \phi_{1}^{2}.
\end{align*}
We know from Proposition \ref{t1} that $\int_{\Omega}a\left(  x\right)
\phi^{2}\log(t_{\ast}\phi_{1})=0$. Thus, we obtain
\begin{equation}
\Phi_{qt}(1,t_{\ast})=\int_{\Omega}a\left(  x\right)  \phi_{1}^{2}>0,
\label{Phqt1t*}%
\end{equation}
as desired. In view of \eqref{Phq1t*}, \eqref{Pht:tt1t*} and \eqref{Phqt1t*},
the Morse Lemma \cite[Theorem 4.3.19]{DM13} applies to $\Phi\in C^{2}$ at
$(1,t_{\ast})$, so that the solution set of $\Phi(q,t)=0$ around $(1,t_{\ast
})$ consists of two continuous curves intersecting transversally at
$(1,t_{\ast})$, as shown in \cite[page 187]{DM13}, where one of them is
$\Gamma_{1}$, and the other one is $\Gamma_{2}$. Since the solution set of
$(P_{a,q})$ is given exactly by $\Gamma_{1}\cup\Gamma_{2}$ around $(1,t_{\ast
})$, assertions (i) and (ii) follow. Finally, assertion (iii) follows from the
fact that $t_{\ast}\phi_{1}\in\mathcal{P}^{\circ}$ and the continuity of
$\Gamma_{2}$ in $\mathbb{R}\times W_{D}^{2,\eta}(\Omega)$ with $\eta>N$.
$\blacksquare$ \newline


\begin{proposition}
\label{prop:gbf} Suppose that $(H_{1})$ holds. Then, $(P_{a,q})$ possesses a
bifurcation curve $\{(q,u(q))\in(\underline{q},1]\times W_{D}^{2,r}(\Omega)\}$
of solutions, for some $\underline{q}\in\lbrack0,1)$, emanating from
$\Gamma_{1}$ at $(1,t_{\ast}\phi_{1})$.
Moreover:

\begin{enumerate}
\item $u(q)\in\mathcal{P}^{\circ}$ for $q\in(\underline{q},1]$, and
$u(1)=t_{\ast}\phi_{1}$;

\item If $\underline{q}>0$ then
\[
\lim_{q\rightarrow\underline{q}^{+}}u(q)=\underline{u}\quad\mbox{in}\ W_{D}%
^{2,r}(\Omega),
\]
and $\underline{u}$ is a nontrivial nonnegative solution of $(P_{a,\underline
{q}})$ which does not belong to $\mathcal{P}^{\circ}$ and satisfies
$\underline{u}>0$ in $\Omega_{+}$.

\item If $\underline{q}=0$ then, for some $\theta\in(0,1)$,
\[
\lim_{q\rightarrow\underline{q}^{+}}u(q)=\underline{u}\not \equiv
0\quad\mbox{ in  }\ C_{0}^{1+\theta}(\overline{\Omega}).
\]
Moreover, $\underline{u}\in W_{\mathrm{loc}}^{2,r}(\Omega_{+})$ satisfies
$\underline{u}>0$ in $\Omega_{+}$.
\end{enumerate}

Furthermore, $(P_{a,q})$ has no other bifurcation points for solutions in
$\mathcal{P}^{\circ}$ on $\Gamma_{1}$, with respect to the $\mathbb{R}\times
W_{D}^{2,\eta}(\Omega)$-topology.
\end{proposition}

To establish Proposition \ref{prop:gbf}, we use the following two lemmas:


\begin{lemma}
[\textit{A priori} estimate]\label{lem:upper} Given $q_{0}\in(0,1)$, there
exists $C_{0}>0$ such that $\Vert u\Vert_{W_{D}^{2,r}(\Omega)}\leq C_{0}$ for
all nontrivial nonnegative solutions of $(P_{a,q})$ with $q\in(0,q_{0}]$.
\end{lemma}

\textit{Proof.} Let $u$ be a nontrivial nonnegative solution of $(P_{a,q})$
with $q\in(0,q_{0}]$. As usual, let $r^{\prime}:=\frac{r}{r-1}>1$. It holds
that
\begin{align*}
\Vert u\Vert_{H_{0}^{1}(\Omega)}^{2}  &  =\int_{\Omega}|\nabla u|^{2}%
=\int_{\Omega}au^{q+1}\\
&  \leq\left\Vert a\right\Vert _{L^{r}(\Omega)}\left\Vert u\right\Vert
_{L^{\left(  q+1\right)  r^{\prime}}(\Omega)}^{q+1}\\
&  \leq C\left\Vert a\right\Vert _{L^{r}(\Omega)}\left\Vert u\right\Vert
_{L^{2^{\ast}}(\Omega)}^{q+1}\\
&  \leq C^{\prime}\left\Vert a\right\Vert _{L^{r}(\Omega)}\left\Vert
u\right\Vert _{H_{0}^{1}(\Omega)}^{q+1}.
\end{align*}
Here,
we have used the fact that $2^{\ast}\geq\frac{(q+1)r}{r-1}$ when $0<q<1$ and
$r>N$. Since $q\leq q_{0}<1$, we deduce that if $\Vert u\Vert_{H_{0}%
^{1}(\Omega)}\geq1$, then
\[
\Vert u\Vert_{H_{0}^{1}(\Omega)}\leq\left(  C^{\prime}\Vert a\Vert
_{L^{r}(\Omega)}\right)  ^{\frac{1}{1-q_{0}}}.
\]
Hence, the boundedness in $H_{0}^{1}(\Omega)$ is verified. An elliptic
regularity argument yields the desired conclusion. $\blacksquare$ \newline


\begin{lemma}
[\textit{A priori} bound from below]\label{lem:zero} Let $q\in\left(
0,1\right)  $ and $\Omega^{\prime}\subseteq\Omega_{+}$ be a smooth domain. If
$u$ is a positive solution of $(P_{a,q})$ then
\begin{equation}
u\geq\lambda_{1}\left(  a,\Omega^{\prime}\right)  ^{-1/\left(  1-q\right)
}\phi\text{\quad in }\Omega^{\prime}, \label{l7}%
\end{equation}
where $\phi>0$ with $\Vert\phi\Vert_{\infty}=1$ is the eigenfunction
associated to $\lambda_{1}(a,\Omega^{\prime})$.

\end{lemma}

\textit{Proof.} Let $\psi$ be given by $\psi:=\lambda_{1}\left(
a,\Omega^{\prime}\right)  ^{-1/\left(  1-q\right)  }\phi$ in $\Omega^{\prime}$
and $\psi:=0$ on $\overline{\Omega}\setminus\Omega^{\prime}$. A few
computations show that $\psi$ is a nonnegative weak subsolution of $\left(
P_{a,q}\right)  $. Since $k\mathcal{S}\left(  a^{+}\right)  $ with
$k\geq\left\Vert \mathcal{S}\left(  a^{+}\right)  \right\Vert _{\infty
}^{1/\left(  1-q\right)  }$ is a supersolution of this problem and $\left(
P_{a,q}\right)  $ has at most one positive solution, we conclude that
$u\geq\psi$ (otherwise, since $\max\left(  u,\psi\right)  $ is a positive weak
subsolution of $\left(  P_{a,q}\right)  $, the method of weak sub and
supersolutions (e.g. \cite[Theorem 4.9]{du}) would yield some $v\not =u$
solution of $\left(  P_{a,q}\right)  $). $\blacksquare$ \newline

\textit{Proof of Proposition \ref{prop:gbf}.} From Proposition \ref{prop:gam2}%
, let $\mathcal{C}_{1}=\{(q,u)\in(0,1]\times W_{D}^{2,\eta}(\Omega)\}$ be the
solution curve of $(P_{a,q})$ bifurcating at $(1,t_{\ast}\phi_{1})$. By
Proposition \ref{prop:open}, we know that $\mathcal{C}_{1}$ consists of
solutions in $\mathcal{P}^{\circ}$ of $(P_{a,q})$ and is parametrized by $q$,
i.e. $u=u(q)$, for $0<q<1$. Moreover, we can define
\[
\underline{q}:=\inf\left\{  q\in
(0,1):\mbox{ $(P_{a,q})$ has a solution $u$ such that $(q,u)
\in \mathcal{C}_1$}\right\}  .
\]
If $\underline{q}>0$, then Proposition \ref{prop:open} does not allow us to
have a solution $\underline{u}$ of $(P_{a,\underline{q}})$ such that
$\underline{u}\in\mathcal{P}^{\circ}$ and $(\underline{q},\underline{u}%
)\in\mathcal{C}_{1}$. For $\mathcal{C}_{1}=\{(q,u(q)):q\in(\underline
{q},1)\}\cup\{(1,t_{\ast}\phi_{1})\}$, Corollary \ref{cor:r} provides that
$q\mapsto u(q)$ is continuous from $(\underline{q},1)$ to $W_{D}^{2,r}%
(\Omega)$.

Now, we investigate the behavior of $u(q)$ as $q\rightarrow1^{-}$. We know
from Proposition \ref{prop:gam2} that $u(q)$ is bounded from both above and
below in $W_{D}^{2,\eta}(\Omega)$ as $q\rightarrow1^{-}$. Since $W_{D}%
^{2,\eta}(\Omega)\subset C_{0}^{1}(\overline{\Omega})$ is compact, if
$q_{n}\rightarrow1^{-}$ then, up to a subsequence, $u(q_{n})\rightarrow w_{1}$
in $C_{0}^{1}(\overline{\Omega})$ for some $w_{1}$, so that $w_{1}=t\phi_{1}$
for some $t>1$. Arguing as in the proof of Corollary \ref{cor:r}, we deduce
that $u(q_{n})\rightarrow t\phi_{1}$ in $W_{D}^{2,r}(\Omega)$. In fact,
Proposition \ref{prop:t*} provides that $t=t_{\ast}$. This argument allows us
to deduce that $u(q)\rightarrow t_{\ast}\phi_{1}$ in $W_{D}^{2,r}(\Omega)$ as
$q\rightarrow1^{-}$. Hence, setting $u(1)=t_{\ast}\phi_{1}$, the map $q\mapsto
u(q)$ is continuous from $(\underline{q},1]$ to $W_{D}^{2,r}(\Omega)$.

Next, we consider the case $\underline{q}>0$ and investigate the behavior of
$u(q)$ as $q\rightarrow\underline{q}^{+}$. Arguing as above, Lemma
\ref{lem:upper} shows that there exists $\underline{u}:=\lim_{q\rightarrow
\underline{q}^{+}}u(q)$ in $C_{0}^{1}(\overline{\Omega})$, using again that
$W_{D}^{2,r}(\Omega)\subset C_{0}^{1}(\overline{\Omega})$ is compact. Thanks
to Lemma \ref{lem:zero}, we have $\underline{u}>0$ in $\Omega_{+}$. By
elliptic regularity, $\underline{u}$ is a nontrivial nonnegative solution of
$(P_{a,\underline{q}})$. As already seen, $\underline{u}\not \in
\mathcal{P}^{\circ}$. Assertion (ii) is now verified.

The case $\underline{q}=0$ can be handled in a similar way as in (ii), using
Lemmas \ref{lem:upper} and \ref{lem:zero}, and therefore (iii) follows.

Finally, the nonexistence result for bifurcation points on $\Gamma_{1}$ is a
direct consequence of Proposition \ref{prop:t*}. $\blacksquare$

\begin{remark}
Let us observe that more information is available for $\underline{u}$ in both
Proposition \ref{prop:gbf} (ii) and (iii). More precisely, let $\Omega_{+}%
^{c}$ be any connected component of $\Omega_{+}$. If $\Omega_{+}^{c}$ is
smooth then, letting $q\rightarrow\underline{q}^{+}$ in (\ref{l7}), it follows
that%
\begin{equation}
\underline{u}\geq\lambda_{1}\left(  a,\Omega_{+}^{c}\right)
^{-1/(1-\underline{q})}\phi_{c}\text{\quad in }\Omega_{+}^{c}, \label{+c}%
\end{equation}
where $\phi_{c}>0$ with $\left\Vert \phi_{c}\right\Vert _{\infty}=1$ is the
positive eigenfunction associated to the weight $a$ in $\Omega_{+}^{c}$. If
$\Omega_{+}^{c}$ is not smooth, then (\ref{+c}) holds with $\Omega_{+}^{c}$
replaced by any smooth domain $\Omega^{\prime}\subset\Omega_{+}$.
\end{remark}

\begin{proposition}
\label{prop:0to1} Under the conditions of Proposition \ref{prop:gbf}, assume
that $\mathcal{S}(a)\in\mathcal{P}^{\circ}$ and $\mathcal{I}_{a}=(0,1)$.
Then $u(q)$ satisfies the following two assertions
(see Figure \ref{fig17_0618b}):


\begin{enumerate}

\item If we set $u(0):=\mathcal{S}(a)$ and $u(1):=t_{\ast}\phi_{1}$, then the
mapping $q\mapsto u(q)$ is continuous from $[0,1]$ to $W_{D}^{2,r}(\Omega)$;

\item $u(q)$ is asymptotically stable for $q\in\lbrack0,1)$.
\end{enumerate}
\end{proposition}

\textit{Proof.}
Corollary \ref{cor:r} provides a curve $\mathcal{C}_{0}=\{(q,u(q))\in
\lbrack0,q_{0})\times W_{D}^{2,r}(\Omega)\}$ of solutions in $\mathcal{P}%
^{\circ}$ of $(P_{a,q})$ emanating from $(0,\mathcal{S}(a))$. On the other
hand, from Proposition \ref{prop:gbf}, we have a curve $\mathcal{C}%
_{1}=\{(q,u(q))\in(q_{1},1]\times W_{D}^{2,r}(\Omega)\}$ of solutions in
$\mathcal{P}^{\circ}$ of $(P_{a,q})$, bifurcating at $(1,t_{\ast}\phi_{1})$.
We set $\overline{q}:=\sup q_{0}$ and $\underline{q}:=\inf q_{1}$. Then,
$0<\overline{q}$ and $\underline{q}<1$. If $\underline{q}<\overline{q}$ then,
by the uniqueness of positive solutions, we see that $\overline{q}=1$,
$\underline{q}=0$ and hence $\mathcal{C}_{0}=\mathcal{C}_{1}$, as desired.

\begin{figure}[tbh]
\begin{center}
\includegraphics[scale=0.25]{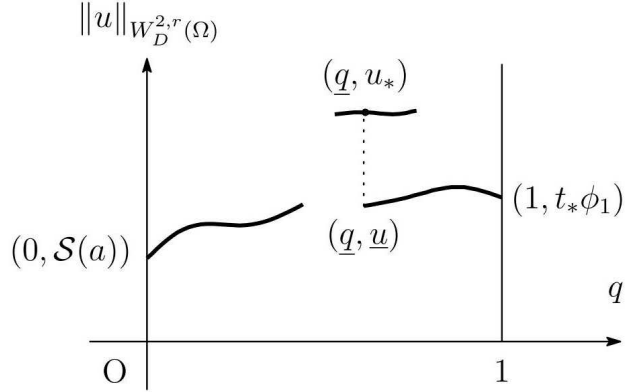}
\end{center}
\caption{The situation of $\underline{u}$ and $u_{*}$.}%
\label{fig17_0618}%
\end{figure}

Assume to the contrary that $\overline{q}\leq\underline{q}$. Then, in
particular, we have $0<\underline{q}<1$. This is the case (ii) in Proposition
\ref{prop:gbf}, according to which we can obtain a nontrivial nonnegative
solution $\underline{u}$ of $(P_{a,\underline{q}})$ such that $\underline
{u}\not \in \mathcal{P}^{\circ}$ and $\underline{u}\in\overline{\mathcal{C}%
_{1}}$. On the other hand, the fact that $\mathcal{I}_{a}=(0,1)$ ensures the
existence of a solution $u_{\ast}\in\mathcal{P}^{\circ}$ of $(P_{a,\underline
{q}})$, which implies $u_{\ast}\neq\underline{u}$. It follows that there exist
$\varepsilon_{\ast}>0$ and some open ball $B_{\ast}$ centered at $u_{\ast}$
such that
\begin{equation}
\left(  (\underline{q},\,\underline{q}+\varepsilon_{\ast})\times B_{\ast
}\right)  \cap\mathcal{C}_{1}|_{q\in(\underline{q},\,\underline{q}%
+\varepsilon_{\ast})}=\emptyset. \label{B*}%
\end{equation}
However, Proposition \ref{prop:open} is also applicable at $(\underline
{q},u_{\ast})$, and then, $(P_{a,q})$ has a solution $u\in\mathcal{P}^{\circ}$
such that $(q,u)\in(\underline{q},\,\underline{q}+\varepsilon_{\ast})\times
B_{\ast}$, which contradicts the uniqueness of positive solutions, see Figure
\ref{fig17_0618}. The assertion (i) is now verified.

Finally, the stability of $u(q)$ is a direct consequence of Proposition
\ref{prop:open}. The proof is complete. $\blacksquare$

Thanks to the homogeneity of $u^{q}$, we obtain the counterpart of Proposition
\ref{prop:0to1} for the case $\lambda_{1}(a)\neq1$ by a change of variables
just as in \cite[Section 3]{neumann}.

\begin{corollary}
\label{cor:0to1} Assume $(H_{1})$, $\mathcal{S}(a)\in\mathcal{P}^{\circ}$ and
$\mathcal{I}_{a}=(0,1)$.
Then the asymptotic profile of $u(q)$
is given by
\[
\lambda_{1}(a)^{\frac{1}{1-q}}u(q)\rightarrow t_{\ast}\phi_{1}\quad
\mbox{in}\quad W_{D}^{2,r}(\Omega),\quad\mbox{as}\quad q\rightarrow1^{-}.
\]
Moreover, the following three assertions hold (see Figure \ref{fig:double}):

\begin{enumerate}
\item Assume that $\lambda_{1}(a)>1$.
If we set $u(0):=\mathcal{S}(a)$ and $u(1):=0$, then the map $q\mapsto u(q)$
is continuous from $[0,1]$ to $W_{D}^{2,r}(\Omega)$;

\item Assume that $\lambda_{1}(a)<1$. Then
\begin{equation}
\Vert u(q)\Vert_{C(\overline{\Omega})}\rightarrow\infty
\ \mbox{ as }q\rightarrow1^{-}. \label{uqinfty}%
\end{equation}
If, in addition, we set $u(0):=\mathcal{S}(a)$, then the map $q\mapsto u(q)$
is continuous from $[0,1)$ to $W_{D}^{2,r}(\Omega)$.

\item $u(q)$ is asymptotically stable for $q\in\lbrack0,1)$.
\end{enumerate}
\end{corollary}

\textit{Proof.} We argue in the same way as in the proof of \cite[Theorem
1.1]{neumann}. It suffices to note that for $q\in(0,1)$, $u$ is a nontrivial
nonnegative solution of $(P_{a,q})$ if and only if $v=\lambda_{1}(a)^{\frac
{1}{1-q}}u$ is a nontrivial nonnegative solution of $(P_{b,q})$ with
$b=\lambda_{1}(a)a$. Since $\lambda_{1}(b)=1$, all the assertions except
\eqref{uqinfty} are straightforward from Proposition \ref{prop:0to1}. For
\eqref{uqinfty}, we deduce from Proposition \ref{prop:0to1} (i)
that $\Vert u(q)\Vert_{W_{D}^{2,r}(\Omega)}\rightarrow\infty$ as
$q\rightarrow1^{-}$. If $\Vert u(q)\Vert_{C(\overline{\Omega})}$ is bounded as
$q\rightarrow1^{-}$, then, by standard elliptic estimates, $\Vert
u(q)\Vert_{W_{D}^{2,r}(\Omega)}$ is bounded as $q\rightarrow1^{-}$, a
contradiction. This yields \eqref{uqinfty}. $\blacksquare$\newline

\begin{remark}
\strut

\begin{enumerate}
\item Proposition \ref{prop:0to1} and Corollary \ref{cor:0to1} provide, for a
weight $a$ such that $\mathcal{S}(a)\in\mathcal{P}^{\circ}$, the whole picture
of the positive solutions set of $(P_{a,q})$ with $q\in(0,1)$.

\item It is easy to find weights $a$ satisfying the assumptions in Proposition
\ref{prop:0to1} and/or Corollary \ref{cor:0to1}. Indeed, let $a_{1},a_{2}\in
C(\overline{\Omega})$ with $a_{1},a_{2}\geq0$ nontrivial, and such that
$a_{1}$ satisfies the decay condition in $(H_{1})$. Then, for all
$\lambda,\varepsilon>0$ small enough, $a_{\lambda,\varepsilon}:=a_{1}-\lambda
a_{2}\chi_{\Omega\setminus\Omega_{\varepsilon}}$ fulfills the hypothesis in
Corollary \ref{cor:0to1}, while a suitable multiple of $a_{\lambda
,\varepsilon}$ satisfies the assumptions of Proposition \ref{prop:0to1}.

\item Assume
that $(P_{a,q})$ has a unique positive solution $u(q)$ for every $q\in(0,1)$,
and $(H_{+})$ holds. Then, the map $q\mapsto u(q)$ is continuous from $(0,1)$
to $W_{D}^{2,r}(\Omega)$. To verify this, given $q_{1}\in(0,1)$, let us
consider the limiting behavior of $u(q)$ as $q\rightarrow q_{1}$. In the same
way as in Proposition \ref{prop:gbf} (ii), any sequence $q_{n}\rightarrow
q_{1}$ has a convergent subsequence, still denoted as $q_{n}$, such that
$u(q_{n})\rightarrow u_{1}$ in $W_{D}^{2,r}(\Omega)$ for some nontrivial
nonnegative solution $u_{1}$ of $(P_{a,q_{1}})$ satisfying that $u_{1}>0$ in
$\Omega_{+}$. Recalling Remark \ref{r0} we conclude that $u_{1}=u(q_{1})$.
This argument yields the continuity of $u(q)$ at $q=q_{1}$, as desired.

From this observation, we find that if we assume, in addition to $(H_{+})$,
the conditions of Proposition \ref{prop:0to1} with $\mathcal{I}_{a}=(0,1)$
replaced by the existence of a positive solution of $(P_{a,q})$ for every
$q\in(0,1)$, then the conclusions of Proposition \ref{prop:0to1} still hold,
with now $u(q)$ asymptotically stable for $q$ close to $0$ and $1$, see Figure
\ref{fig17_0623a-m2}.

\begin{figure}[tbh]
\begin{center}
\includegraphics[scale=0.25]{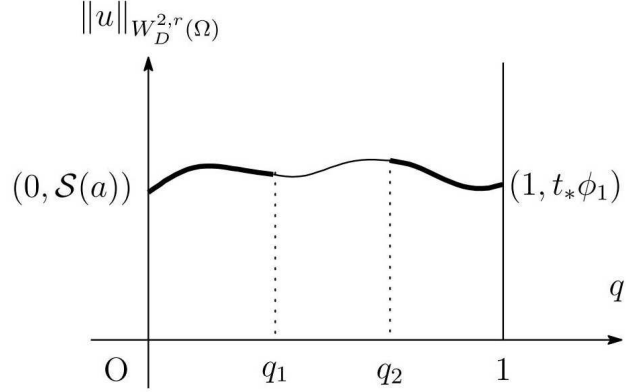}
\end{center}
\caption{The curve of positive solutions in the case $\lambda_{1}(a)=1$.}%
\label{fig17_0623a-m2}%
\end{figure}
\end{enumerate}
\end{remark}

\section{Further results and proofs of the main theorems}

\label{sec:p}

We set
\begin{align*}
\mathcal{A}_{a}  &  :=\{q\in(0,1):\text{$(P_{a,q})$ has a unique nontrivial
nonnegative solution,}\\
&  \text{which in addition belongs to }\mathcal{P}^{\circ}\},
\end{align*}
and recall that%
\[
\mathcal{I}_{a}:=\{q\in(0,1):(P_{a,q})\text{ has a solution }u\in
\mathcal{P}^{\circ}\}.
\]
As a consequence of \cite[Corollary 1.5 and Theorem 1.9]{krqu} we know that
$\mathcal{A}_{a}=\left(  q_{0},1\right)  $ for some $q_{0}\geq0$. The next
proposition shows that in general $\mathcal{I}_{a}\not =\mathcal{A}_{a}$ if
the first condition in $(H_{+}^{\prime})$ is not satisfied. It also shows that
the ground state solution $U_{q}$ is not necessarily a positive solution of
$(P_{a,q})$.

\begin{proposition}
\label{propo}Let $\Omega:=(x_{0},x_{1})\subset\mathbb{R}$ and $q\in\left(
0,1\right)  $.\strut

\begin{enumerate}
\item There exists $a\in C(\Omega)\cap L^{r}(\Omega)$, $r>1$, such that
$q\in\mathcal{I}_{a}\setminus\mathcal{A}_{a}$.

\item There exists $a\in C(\Omega)\cap L^{r}(\Omega)$, $r>1$, such that the
nonnegative ground state solution of $(P_{a,q})$ vanishes in a subinterval of
$\Omega$.

\end{enumerate}
\end{proposition}

\textit{Proof}. After a translation and a dilation, we may assume that
$\Omega:=\left(  -2,2\right)  $.

\begin{figure}[tbh]
\begin{center}
\includegraphics[scale=0.25]{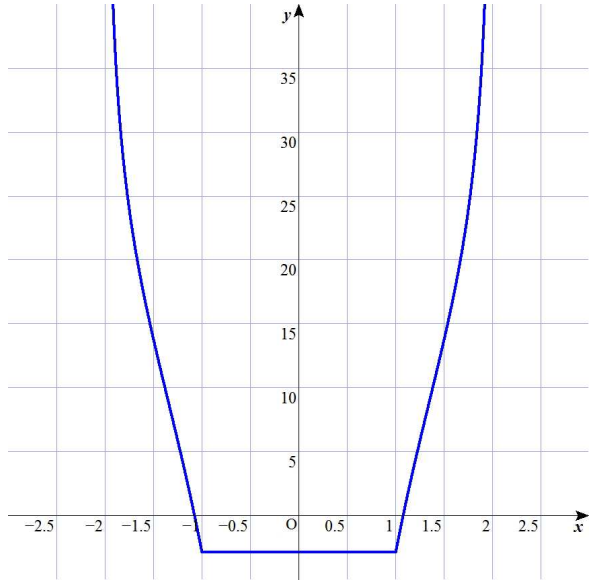}
\end{center}
\caption{The indefinite weight $a$ in the case $q= \frac{1}{3}$.}%
\label{fig17_0604b}%
\end{figure}

\begin{figure}[tbh]
\begin{center}
\includegraphics[scale=0.35]{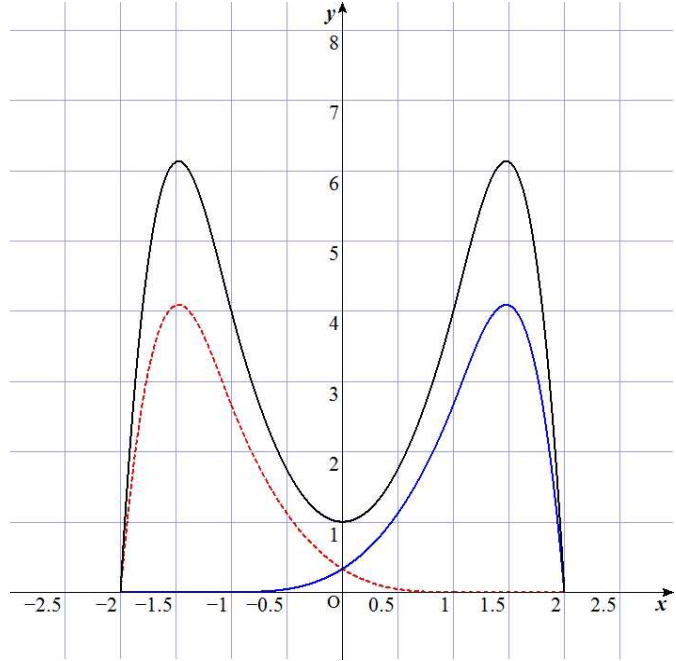}
\end{center}
\caption{The nontrivial nonnegative solutions $u_{1},u_{2}$ with dead cores
and a possible solution $u_{3}$ in $\mathcal{P}^{\circ}$ (which is even, by
uniqueness) in the case $q=\frac{1}{3}$.}%
\label{fig17_0604a}%
\end{figure}

\begin{enumerate}
\item Let $q\in\left(  0,1\right)  $. We first define
\[
r:=\frac{2}{1-q}\in\left(  2,\infty\right)  \quad\text{and}\quad f\left(
x\right)  :=\frac{\left(  x+1\right)  ^{r}}{r}.
\]
Note that $rq=r-2$. Let $p$ be the polynomial given by
\[
p\left(  x\right)  :=\alpha x^{3}+\beta x^{2}+\gamma x+\delta,
\]
where
\begin{gather*}
\alpha:=-2^{r-3}\left(  \frac{8}{r}+r+3\right)  ,\quad\beta:=2^{r-1}\left(
r+\frac{6}{r}+2\right)  ,\\
\gamma:=-2^{r-3}\left(  5r+\frac{24}{r}+3\right)  ,\quad\delta:=2^{r-2}\left(
\frac{8}{r}+r-1\right)  .
\end{gather*}
One can verify that
\begin{equation}
p\left(  1\right)  =f\left(  1\right)  ,\quad p^{\prime}\left(  1\right)
=f^{\prime}\left(  1\right)  ,\quad p^{\prime\prime}\left(  1\right)
=f^{\prime\prime}\left(  1\right)  ,\quad p\left(  2\right)  =0. \label{mm}%
\end{equation}
Moreover, $p^{\prime}\left(  2\right)  =-2^{r-3}\left(  7+r+\frac{24}%
{r}\right)  <0$. Thus, since $p\left(  1\right)  ,p^{\prime}\left(  1\right)
>0$, $p$ has degree $3$ and $\alpha<0$, it follows that $p>0$ in $\left(
1,2\right)  $ (because, if not, $p^{\prime}$ would vanish at least at three
different points, which is not possible).

We now set%
\[
a\left(  x\right)  :=%
\begin{cases}
-\left(  r-1\right)  r^{q} & \text{if }x\in\left[  0,1\right]  ,\\
-\frac{p^{\prime\prime}\left(  x\right)  }{\left[  p\left(  x\right)  \right]
^{q}} & \text{if }x\in\left[  1,2\right)  ,\\
a\left(  -x\right)  & \text{if }x\in\left(  -2,0\right]  .
\end{cases}
\]
Observe that $a\in C(\Omega)$ since $rq=r-2$ and so
\[
-\frac{p^{\prime\prime}\left(  1\right)  }{\left[  p\left(  1\right)  \right]
^{q}}=-\frac{f^{\prime\prime}\left(  1\right)  }{\left[  f\left(  1\right)
\right]  ^{q}}=-\left(  r-1\right)  r^{q}.
\]
Also, since $p\left(  x\right)  ^{q}$ behaves like $C\left(  x-2\right)  ^{q}$
as $x\rightarrow2^{-}$ we derive that $a\in L^{r}(\Omega)$ for some $r>1$, see
Figure \ref{fig17_0604b}.

Define now
\[
u_{1}\left(  x\right)  :=%
\begin{cases}
0 & \text{if }x\in\left[  -2,-1\right]  ,\\
f\left(  x\right)  & \text{if }x\in\left[  -1,1\right]  ,\\
p\left(  x\right)  & \text{if }x\in\left[  1,2\right]  ,
\end{cases}
\]
and $u_{2}\left(  x\right)  :=u_{1}\left(  -x\right)  $. By (\ref{mm}) we have
that $u_{1},u_{2}\in C^{2}(\overline{\Omega})$, see Figure \ref{fig17_0604a}.
Moreover, one can check that $u_{1}$ and $u_{2}$ are two nonnegative
nontrivial solutions of the problem%
\begin{equation}%
\begin{cases}
-u^{\prime\prime}=a(x)u^{q} & \mbox{ in }\Omega,\\
u=0 & \mbox{ on }\partial\Omega.
\end{cases}
\label{eje}%
\end{equation}

Furthermore, integrating by parts we deduce that
\[
z\left(  x\right)  :=\max\left(  u_{1}\left(  x\right)  ,u_{2}\left(
x\right)  \right)  =\left\{
\begin{array}
[c]{lll}%
u_{2}(x) & \mathrm{in} & \left[  -2,0\right]  ,\\
u_{1}(x) & \mathrm{in} & \left[  0,2\right]  ,
\end{array}
\right.
\]
is a strictly positive weak subsolution of (\ref{eje}), satisfying that
$\left\vert z^{\prime}\left(  \pm2\right)  \right\vert >0$. Thus, since
$k\mathcal{S}\left(  a^{+}\right)  $ is a supersolution of (\ref{eje}) for
every $k>0$ large enough, we obtain then a solution $u\in\mathcal{P}^{\circ}$
of (\ref{eje}). In other words, $q\in\mathcal{I}_{a}$ but $q\not \in
\mathcal{A}_{a}$, since $u_{1}$ and $u_{2}$ are nontrivial nonnegative
solutions vanishing in nonempty subdomains of $\Omega$.

\item We set $\underline{a}\left(  x\right)  :=a\left(  x\right)  $ for
$x\in\left[  -1,2\right]  $ and $\underline{a}\left(  x\right)  :=a\left(
-1\right)  $ for $x\in\left[  -2,-1\right]  $. It is clear that $u_{1}$ is a
solution of $(P_{\underline{a},q})$. Moreover, since $u>0$ in $\Omega
_{+}\left(  \underline{a}\right)  $, it follows from the uniqueness assertion
in Remark \ref{r0} that $u_{1}$ is the ground state solution of
$(P_{\underline{a},q})$. $\blacksquare$
\end{enumerate}

The above proposition also shows that $\mathcal{A}_{a}$ can have arbitrarily
small size, since $\mathcal{A}_{a}=\left(  q_{0},1\right)  $ for some
$q_{0}\geq0$. On the other hand, it was asked in \cite{krqu} whether it could
happen that $\mathcal{A}_{a}=(0,1)$. The next proposition gives a positive
answer to the aforementioned question.

\begin{proposition}
\label{ti}Let $a\in L^{r}\left(  \Omega\right)  $, $r>N$.\strut

\begin{enumerate}
\item Suppose $\mathcal{S}\left(  a\right)  >0$ in $\Omega$, and either
$a\geq0$ in $\Omega_{\rho_{0}}$ for some $\rho_{0}>0$, or $\Omega$ is a ball
and $a$ is radial. Then $\mathcal{I}_{a}=\left(  0,1\right)  $.

\item Assume $\mathcal{S}\left(  a\right)  \in\mathcal{P}^{\circ}$ and
$(H_{+}^{\prime})$ holds. Then $\mathcal{I}_{a}=\mathcal{A}_{a}=\left(
0,1\right)  $.
\end{enumerate}
\end{proposition}

\textit{Proof}.

\begin{enumerate}
\item Let $\phi:=\mathcal{S}\left(  a\right)  $, and $\alpha:=1/\left(
1-q\right)  >1$ for $q\in\left(  0,1\right)  $. Note that $\alpha-1=\alpha q$.
Let $0\leq\zeta\in C_{c}^{\infty}\left(  \Omega\right)  $ and $\Omega^{\prime
}$ be a smooth open set with \textit{supp} $\zeta\subset\Omega^{\prime}%
\Subset\Omega$. It holds that $\phi^{\alpha}\in W_{0}^{1,r}\left(
\Omega\right)  \cap W^{2,r}\left(  \Omega^{\prime}\right)  $ and
\begin{align*}
-\Delta(\frac{\phi}{\alpha})^{\alpha}  &  =-\frac{1}{\alpha^{\alpha-1}}\left(
\phi^{\alpha-1}\Delta\phi+\left(  \alpha-1\right)  \phi^{\alpha-2}\left\vert
\nabla\phi\right\vert ^{2}\right) \\
&  \leq-(\frac{\phi}{\alpha})^{\alpha-1}\Delta\phi=a\left(  x\right)
(\frac{\phi}{\alpha})^{\alpha q}\quad\text{in }\Omega^{\prime}\text{.}%
\end{align*}
Let $\psi:=\phi/\alpha$. Multiplying the above inequality by $\zeta$,
integrating over $\Omega^{\prime}$ and using the divergence theorem we derive
that
\[
\int_{\Omega}\nabla(\psi^{\alpha})\nabla\zeta\leq\int_{\Omega}a\left(
x\right)  \psi^{\alpha q}\zeta.
\]
Now, let $0\leq v\in W_{0}^{1,2}\left(  \Omega\right)  $. There exists
$\left\{  \zeta_{n}\right\}  \subset C_{c}^{\infty}\left(  \Omega\right)  $
with $\zeta_{n}\geq0$ in $\Omega$ and such that $\zeta_{n}\rightarrow v$ in
$W^{1,2}\left(  \Omega\right)  $ (e.g. \cite{chipot}, p. 50). Employing the
above inequality with $\zeta$ replaced by $\zeta_{n}$ and going to the limit,
we infer that $-\Delta(\psi^{\alpha})\leq a\left(  x\right)  \left(
\psi^{\alpha}\right)  ^{q}$ in the weak sense in $\Omega$.

On the other side, it is easy to see that $k\mathcal{S}\left(  a^{+}\right)  $
is a supersolution of $(P_{a,q})$ for all $k>0$ large enough. Thus, the method
of weak sub and supersolutions provides us a weak solution $u$ of $(P_{a,q})$.
Moreover, $u>0$ in $\Omega$ and, since $\mathcal{S}\left(  a^{+}\right)  $ is
bounded, by standard regularity arguments we deduce that $u\in W_{D}%
^{2,r}\left(  \Omega\right)  $.

Now, if $a\geq0$ in $\Omega_{\rho_{0}}$, Hopf's Lemma (e.g. \cite[Proposition
1.16]{jai}) tells us that $\partial u/\partial\nu<0$ on $\partial\Omega$. In
other words, $u\in\mathcal{P}^{\circ}$. Suppose now that $\Omega$ is a ball
and $a$ is radial. By the first part of the proof, $(P_{a,q})$ has a positive
solution $u$. Assume by contradiction that $u\not \in \mathcal{P}^{\circ}$.
Since the positive solution of $(P_{a,q})$ is unique and $a$ is radial, we
have that $u$ is radial. Taking this into account we derive that
$\frac{\partial u}{\partial\nu}=0$ on $\partial\Omega$. In other words, $u$ is
a positive solution of the Neumann problem%
\[%
\begin{cases}
-\Delta u=a(x)u^{q} & \mbox{ in }\Omega,\\
\partial u/\partial\nu=0 & \mbox{ on }\partial\Omega.
\end{cases}
\]
Therefore we have that $\int_{\Omega}a<0$ (see \cite{BPT2} or \cite{neumann}).
But this is not possible, because $\mathcal{S}(a)>0$ implies that
$\int_{\Omega}a\geq0$. Hence $u\in\mathcal{P}^{\circ}$, and the proof of (i)
is concluded.

\item We start showing that $\mathcal{S}\left(  a\right)  \in\mathcal{P}%
^{\circ}$ implies that $\left(  0,q_{0}\right)  \subset\mathcal{I}_{a}$ for
some $q_{0}\in\left(  0,1\right)  $. For $\delta>0$, define%
\[
a_{\delta}\left(  x\right)  :=a^{+}\left(  x\right)  \chi_{\Omega
\setminus\Omega_{\delta}}-\left(  1+\delta\right)  a^{-}\left(  x\right)  .
\]
It holds that $a_{\delta}\rightarrow a$ in $L^{r}\left(  \Omega\right)  $ as
$\delta\rightarrow0^{+}$. Hence, since $\mathcal{S}:L^{r}\left(
\Omega\right)  \rightarrow C^{1}(\overline{\Omega})$ is a continuous operator
and $\mathcal{S}\left(  a\right)  \in\mathcal{P}^{\circ}$, we may fix
$\delta>0$ such that $\phi_{\delta}:=\mathcal{S}\left(  a_{\delta}\right)
\in\mathcal{P}^{\circ}$.

Now, we can also fix a constant $k$ such that
\begin{equation}
\frac{1}{1+\delta}<k<1. \label{k}%
\end{equation}
Taking into account (\ref{k}), we see that there exists $q_{0}>0$ such that
for all $q\in\left(  0,q_{0}\right)  $,%
\begin{align*}
\left[  \frac{\left(  \phi_{\delta}\right)  ^{q}}{1+\delta}\right]
^{1/\left(  1-q\right)  }  &  \leq\left[  \frac{\left\Vert \phi_{\delta
}\right\Vert _{\infty}^{q}}{1+\delta}\right]  ^{1/\left(  1-q\right)  }\leq
k\quad\text{in }\Omega,\\
k  &  \leq\left(  \inf_{\Omega\setminus\Omega_{\delta}}\phi_{\delta}\right)
^{q/\left(  1-q\right)  }\leq\left(  \phi_{\delta}\right)  ^{q/\left(
1-q\right)  }\quad\text{in }\Omega\setminus\Omega_{\delta}\text{.}%
\end{align*}
Therefore, for all such $q$'s,%
\begin{align*}
-\Delta\left(  k\phi_{\delta}\right)   &  =ka^{+}\left(  x\right)
\chi_{\Omega\setminus\Omega_{\delta}}-k\left(  1+\delta\right)  a^{-}\left(
x\right) \\
&  \leq a^{+}\left(  x\right)  \left(  k\phi_{\delta}\right)  ^{q}\chi
_{\Omega\setminus\Omega_{\delta}}-a^{-}\left(  x\right)  \left(  k\phi
_{\delta}\right)  ^{q}\leq a\left(  x\right)  \left(  k\phi_{\delta}\right)
^{q}\quad\text{in }\Omega\text{.}%
\end{align*}
In other words, $k\phi_{\delta}\in\mathcal{P}^{\circ}$ is a subsolution of
$(P_{a,q})$. Thus, arguing as above, it follows that there exists a solution
of $(P_{a,q})$ lying in $\mathcal{P}^{\circ}$ for all $q\in\left(
0,q_{0}\right)  $. So, $\left(  0,q_{0}\right)  \subset\mathcal{I}_{a}$.

Next we observe that, thanks to $(H_{+}^{\prime})$, we have that
$\mathcal{I}_{a}=\mathcal{A}_{a}$. The proof follows the lines of the Neumann
case (see \cite[Theorem 1.8]{neumann}), but we repeat it here for the sake of
completeness. It suffices to see that for $q\in\mathcal{I}_{a}$, there exists
a unique nontrivial nonnegative solution of $(P_{a,q})$. Assume by
contradiction that $u$ and $v$ are nontrivial nonnegative solutions of
$(P_{a,q})$ with $u\not \equiv v$. We note that $u\not \equiv 0$ in
$\Omega_{+}$. Indeed, if not, then $\Delta u\geq0$ in $\Omega$ and $u=0$ on
$\partial\Omega$, and therefore the maximum principle says that $u\equiv0$ in
$\Omega$, which is not possible. Now, since $\Omega_{+}$ is connected by
$(H_{+}^{\prime})$, arguing as in Lemma 2.1 in \cite{BPT} we derive that $u>0$
in $\Omega_{+}$. Moreover, the same reasoning applies to $v$, so that $v>0$ in
$\Omega_{+}$. But on the other side, by the uniqueness assertion in Remark
\ref{r0} there exists at most one solution which is positive in $\Omega_{+}$.
So we reach a contradiction.

To end the proof we note that by \cite[Theorems 1.3 and 1.9]{krqu},
$\mathcal{A}_{a}=(q_{a},1)$ for some $q_{a}\in\left[  0,1\right)  $. So, since
$\left(  0,q_{0}\right)  \subset\mathcal{I}_{a}=\mathcal{A}_{a}$, the proof is
complete. $\blacksquare$
\end{enumerate}

\begin{remark}
\label{rem-ti}\strut

\begin{enumerate}
\item It is an interesting open question whether the set $\mathcal{I}_{a}$ is
connected or not.

\item When $a\in C(\overline{\Omega})$, one can see, by the strong maximum
principle, that the condition $\mathcal{S}\left(  a\right)  >0$ in $\Omega$ is
equivalent to $\mathcal{S}\left(  a\right)  >0$ in $\left\{  x\in
\Omega:a\left(  x\right)  \leq0\right\}  $.

\item Theorem 4.4 in \cite{jesusultimo} says that $(P_{a,q})$ admits a
positive solution $u$ (which, a priori, may not belong to $\mathcal{P}^{\circ
}$) for any $q\in\left(  0,1\right)  $, provided that $\mathcal{S}\left(
a\right)  \in\mathcal{P}^{\circ}$. Since the positive solution of $(P_{a,q})$
is unique, we see that $u$ lies in fact in $\mathcal{P}^{\circ}$ whenever some
of the conditions in Proposition \ref{ti}
are satisfied.

\item One can still obtain that $\mathcal{I}_{a}=\left(  0,1\right)  $ in
Proposition \ref{ti} (ii) without assuming $(H_{+}^{\prime})$, provided that
$a$ is \textquotedblleft sufficiently\textquotedblright\ positive in some
component of $\Omega_{+}$. More precisely, let $\Omega_{+}^{c}$ be any
connected component of $\Omega_{+}$, and let $\Omega^{\prime}\subseteq
\Omega_{+}^{c}$ be a smooth domain. Set $\underline{a}:=0$ in $\Omega
_{+}\setminus\Omega^{\prime}$ and $\underline{a}:=a$ in $\Omega\setminus
\left(  \Omega_{+}\setminus\Omega^{\prime}\right)  $. If $\mathcal{S}\left(
\underline{a}\right)  \in\mathcal{P}^{\circ}$, then Proposition \ref{ti} (ii)
tells us that $\mathcal{I}_{\underline{a}}=\left(  0,1\right)  $. Now, since
$\underline{a}\leq a$, for any $q\in\left(  0,1\right)  $ the positive
solution of $(P_{\underline{a},q})$ is a subsolution of $(P_{a,q})$, and this
in turn implies that $\mathcal{I}_{a}=\left(  0,1\right)  $.

\item In the one-dimensional case, \cite[Proposition 2.5 (i) and Remark 2.6
(i)]{nodea} provide us with a class of $a$ such that $\int_{\Omega}a>0$ and
$(P_{a,q})$ has a positive solution $u\in\mathcal{P}^{\circ}$ for all
$q\in\left(  0,1\right)  $, i.e. $\mathcal{I}_{a}=\left(  0,1\right)  $.
Moreover, it is easy to check that for such $a$ we have $\mathcal{I}%
_{a}=\mathcal{A}_{a}$. Finally, one can choose, in this class, some $a$ such
that $\mathcal{S}\left(  a\right)  \not \in \mathcal{P}^{\circ}$. It follows
that the condition $\mathcal{S}\left(  a\right)  \in\mathcal{P}^{\circ}$ is
\textit{not }necessary in order to have $\mathcal{A}_{a}=\left(  0,1\right)
$. Note that, by the divergence theorem, if $\mathcal{S}\left(  a\right)
\in\mathcal{P}^{\circ}$ then $\int_{\Omega}a>0$, i.e., the condition
$\mathcal{S}\left(  a\right)  \in\mathcal{P}^{\circ}$ is stronger than
$\int_{\Omega}a>0$. Let us observe, however, that $\int_{\Omega}a>0$ is
\textit{neither} necessary nor sufficient to have $\mathcal{A}_{a}=\left(
0,1\right)  $. Indeed, Proposition \ref{propo} shows that it is not sufficient
(see also Figure \ref{fig17_0604b}). On the other hand, define
\[
a\left(  x\right)  :=1-2\cos^{2}(x-\frac{\pi}{2}),\quad x\in\Omega
:=(-\frac{\pi}{2},\frac{\pi}{2}).
\]
Then, $a$ is even and $\mathcal{S}\left(  a\right)  =\frac{1}{2}\sin
^{2}(x-\frac{\pi}{2})>0$ in $\Omega$. Thus, Proposition \ref{ti} (i) says that
$\mathcal{I}_{a}=\left(  0,1\right)  $. Also, since $(H_{+}^{\prime})$ is
satisfied, the proof of Proposition \ref{ti} (ii) gives that $\mathcal{I}%
_{a}=\mathcal{A}_{a}$. We note, however, that $\int_{\Omega}a=0$. Note also
that $\mathcal{S}(a)\not \in \mathcal{P}^{\circ}$.
\end{enumerate}
\end{remark}

\quad

We conclude the paper providing the proofs of our main theorems.

\quad

\textit{Proof of Theorem \ref{mt2a}. }The first assertions and (i) are a
direct consequence of Propositions \ref{pgs} and \ref{pgs2}. Item (ii) follows
from Propositions \ref{prop:gam2} and \ref{prop:gbf}, while the first
assertion in (iii) follows as in \cite[Theorem 2.3 (ii)]{BPT}. Moreover, by
Remark \ref{remgs} (i) we have \eqref{iia}. Finally, the fact that $I_{a}$ is
open follows from the continuity of $q\mapsto U_{q}$ from $(0,1)$ to
$C_{0}^{1}(\overline{\Omega})$. $\blacksquare$

\quad

\textit{Proof of Theorem \ref{mt2b}. }The existence of a (unique) positive
solution $u(q)$ is a consequence of the first part of the proof of Proposition
\ref{ti}. Moreover, the aforementioned proof shows that $\left[  \left(
1-q\right)  \mathcal{S}\left(  a\right)  \right]  ^{1/\left(  1-q\right)  }$
is a positive subsolution for $\left(  P_{a,q}\right)  $ and so we get that
\[
u\left(  q\right)  \geq\left[  \left(  1-q\right)  \mathcal{S}\left(
a\right)  \right]  ^{1/\left(  1-q\right)  }>0\text{\quad in }\Omega
\]
(otherwise, arguing as in the proof of Lemma \ref{lem:zero}, we would
contradict the uniqueness of Theorem 1.0 (i)). Now, taking into account this
inequality and the \textit{a priori} bound in Lemma \ref{lem:zero} we derive
that $u(q)\rightarrow\mathcal{S}(a)$ in $W_{D}^{2,r}(\Omega)$ as
$q\rightarrow0^{+}$.

Next, we note that item (i) is a consequence of Theorem \ref{mt2a} and the
uniqueness assertion in Remark \ref{r0}. On the other side, (ii) follows from
Theorem \ref{mt2a} (i) and the first part of the proof of Proposition \ref{ti}
(ii), while we may conclude (iii) from the aforementioned proposition.
$\blacksquare$

\quad

\textit{Proofs of Theorem \ref{mt2c} and Corollary \ref{sing}. }Both are
immediate from Proposition \ref{prop:open}. $\blacksquare$

\section*{Acknowledgements}

The authors thank the referee for carefully reading the paper and pointing out
some corrections. U. Kaufmann was partially supported by Secyt-UNC
33620180100016CB; H. Ramos Quoirin was supported by Fondecyt grants 1161635,
1171532, 1171691, and 1181125; K. Umezu was supported by JSPS KAKENHI Grant Numbers JP15K04945 and JP18K03353.


\begin{thebibliography}{99}                                                                                               %


\bibitem {AL09}S.\ Alama, Q.\ Lu, \textit{Compactly supported solutions to
stationary degenerate diffusion equations,} J.\ Differential Equations
\textbf{246} (2009), 3214--3240.

\bibitem {BPT}C.\ Bandle, M.\ Pozio, A.\ Tesei, \textit{The asymptotic
behavior of the solutions of degenerate parabolic equations},
Trans.\ Amer.\ Math.\ Soc.\ \textbf{303} (1987), 487--501.


\bibitem {BPT2}C.\ Bandle, A.\ M.\ Pozio, A.\ Tesei, \textit{Existence and
uniqueness of solutions of nonlinear Neumann problems},
Math.\ Z.\ \textbf{199} (1988), 257--278.

\bibitem {BK92}H.\ Brezis, S.\ Kamin, \textit{Sublinear elliptic equations in
}$\mathbb{R}^{n}$, Manuscripta Math.\ \textbf{74} (1992), 87--106.

\bibitem {chipot}M.\ Chipot, \textit{Elliptic equations: an introductory
course}. Birkh\"{a}user Advanced Texts, Birkh\"{a}user Verlag, Basel, 2009.

\bibitem {tartar}M.\ Crandall, P.\ Rabinowitz, L.\ Tartar, \textit{On a
Dirichlet problem with a singular nonlinearity}, Comm.\ Partial Differential
Equations \textbf{2} (1977), 193--222.

\bibitem {cuesta}M.\ Cuesta, P.\ Tak\'{a}\v{c}, \textit{A strong comparison
principle for positive solutions of degenerate elliptic equations},
Differential Integral Equations \textbf{13} (2000), 721--746.

\bibitem {jai}D.\ De Figueiredo, \textit{Positive solutions of semilinear
elliptic equations}, Lect.\ Notes Math.\ Springer \textbf{957} (1982), 34--87.

\bibitem {DS}M.\ Delgado, A.\ Su\'{a}rez, \textit{On the uniqueness of
positive solution of an elliptic equation}, Appl. Math. Lett. \textbf{18}
(2005), 1089--1093.

\bibitem {pino}M.\ del Pino, \textit{A global estimate for the gradient in a
singular elliptic boundary value problem}, Proc.\ Roy.\ Soc.\ Edinburgh
Sect.\ A \textbf{122} (1992), 341--352.

\bibitem {DM13}P.\ Dr\'{a}bek, J.\ Milota, Methods of nonlinear analysis,
Applications to differential equations, Second edition, Birkh\"{a}user
Advanced Texts: Basler Lehrb\"{u}cher, Birkh\"{a}user/ Springer Basel AG,
Basel, 2013.

\bibitem {du}Y. Du, \textit{Order structure and topological methods in
nonlinear partial differential equations. Vol.\ 1. Maximum principles and
applications}, World Scientific Publishing Co. Pte. Ltd., Hackensack, NJ, 2006.

\bibitem {jesusultimo}J. Hern\'{a}ndez, F. Mancebo, J. Vega, \textit{On the
linearization of some singular, nonlinear elliptic problems and applications},
Ann. Inst. H. Poincar\'{e} Anal. Non Lin\'{e}aire \textbf{19} (2002), 777--813.

\bibitem {royal}J. Hern\'{a}ndez, F. Mancebo, J. Vega, \textit{Positive
solutions for singular nonlinear elliptic equations}, Proc. Roy. Soc.
Edinburgh Sect. A \textbf{137 }(2007), 41--62.

\bibitem {nodea}T. Godoy, U. Kaufmann, \textit{On strictly positive solutions
for some semilinear elliptic problems}, NoDEA Nonlinear Differ. Equ. Appl.
\textbf{20} (2013), 779--795.

\bibitem {jmaa}T. Godoy, U. Kaufmann, \textit{On Dirichlet problems with
singular nonlinearity of indefinite sign}, J. Math. Anal. Appl. \textbf{428}
(2015), 1239--1251.

\bibitem {ans}T. Godoy, U. Kaufmann, \textit{Existence of strictly positive
solutions for sublinear elliptic problems in bounded domains}, Adv. Nonlinear
Stud. \textbf{14} (2014), 353--359.

\bibitem {gomes}S. Gomes, \textit{On a singular nonlinear elliptic problem},
SIAM J. Math. Anal. \textbf{17} (1986), 1359--1369.

\bibitem {J}L. Jeanjean, \textit{Some continuation properties via minimax
arguments}, Electron. J. Differential Equations \textbf{2011} (2011), Paper
No. 48, 10 pp.

\bibitem {K}R. Kajikiya, \textit{Positive solutions of semilinear elliptic
equations with small perturbations}, Proc. Amer. Math. Soc. \textbf{141}
(2013), 1335--1342.

\bibitem {ana}U. Kaufmann, I. Medri, \textit{One-dimensional singular problems
involving the }$p$\textit{-Laplacian and nonlinearities indefinite in sign},
Adv. Nonlinear Anal. \textbf{5} (2016), 251--259.

\bibitem {krqu}U. Kaufmann, H. Ramos Quoirin, K. Umezu,\ \textit{Positivity
results for indefinite sublinear elliptic problems via a continuity argument},
J. Differential Equations \textbf{263} (2017), 4481--4502.


\bibitem {neumann}U. Kaufmann, H. Ramos Quoirin, K. Umezu,\ \textit{Positive
solutions of an elliptic Neumann problem with a sublinear indefinite
nonlinearity. } Preprint. (arXiv:1705.07791)

\bibitem {KLP}B. Kawohl, M. Lucia, S. Prashanth, \textit{Simplicity of the
principal eigenvalue for indefinite quasilinear problems}. Adv. Differential
Equations \textbf{12} (2007), 407--434.

\bibitem {lair}A. Lair, A. Shaker, \textit{Classical and weak solutions of a
singular semilinear elliptic problem}, J. Math. Anal. Appl. \textbf{211}
(1997), 371--385.

\bibitem {lazer}A. Lazer, P. McKenna, \textit{On a singular nonlinear elliptic
boundary-value problem}, Proc. Amer. Math. Soc. \textbf{111 }(1991),
721--730.


\bibitem {PT}M. A. Pozio and A. Tesei, \textit{Support properties of solution
for a class of degenerate parabolic problems}, Comm. Partial Differential
Equations \textbf{12} (1987), 47--75.

\bibitem {Sp83}J.\ Spruck, \textit{Uniqueness in a diffusion model of
population biology}, Comm.\ Partial Differential Equations \textbf{8} (1983), 1605--1620.

\bibitem {Zei86}E. Zeidler, Nonlinear functional analysis and its applications
I: Fixed-point theorems, Translated from the German by Peter R. Wadsack,
Springer-Verlag, New York, 1986.
\end{thebibliography}
\end{document}